\documentclass [10pt,reqno]{amsart}
\usepackage {amsmath,amssymb,verbatim,geometry}
\usepackage{mathtools,esint}
\usepackage[pdftex,hyperindex]{hyperref}

\def\XXint#1#2#3{{\setbox0=\hbox{$#1{#2#3}{\int}$ }
\vcenter{\hbox{$#2#3$ }}\kern-.6\wd0}}

\newcommand{\C}{\mathbb{C}}
\newcommand{\DD}{\mathbb{D}}

\newcommand{\N}{\mathbb{N}}
\renewcommand{\P}{\mathbb{P}}
\newcommand{\Q}{\mathbb{Q}}
\newcommand{\R}{\mathbb{R}}
\newcommand{\Z}{\mathbb{Z}}

\newcommand{\fa}{\mathfrak{a}}
\newcommand{\fb}{\mathfrak{b}}
\newcommand{\fm}{\mathfrak{m}}

\newcommand{\tX}{\widetilde{X}}

\newcommand{\cE}{\mathcal{E}}
\newcommand{\cF}{\mathcal{F}}

\newcommand{\cH}{\mathcal{H}}

\newcommand{\cJ}{\mathcal{J}}

\newcommand{\cL}{\mathcal{L}}
\newcommand{\cM}{\mathcal{M}}

\newcommand{\cO}{\mathcal{O}}

\newcommand{\cX}{\mathcal{X}}

\newcommand{\Xan}{X^{\mathrm{an}}}

\newcommand{\Xdiv}{X^{\mathrm{div}}}

\renewcommand{\a}{\alpha}
\renewcommand{\b}{\beta}
\renewcommand{\d}{\delta}
\newcommand{\e}{\varepsilon}
\newcommand{\f}{\varphi}

\newcommand{\la}{\lambda}
\newcommand{\om}{\omega}
\newcommand{\p}{\psi}

\newcommand{\eg}{{\rm e.g.\ }} 
\newcommand{\ie}{{\rm i.e.\ }} 
\newcommand{\triv}{\mathrm{triv}} 

\newcommand{\loc}{\mathrm{loc}}

\DeclareMathOperator{\opD}{D}
\DeclareMathOperator{\opE}{E}

\DeclareMathOperator{\opF}{F}
\DeclareMathOperator{\opH}{H}
\DeclareMathOperator{\opI}{I}
\DeclareMathOperator{\opJ}{J}
\DeclareMathOperator{\opL}{L}
\DeclareMathOperator{\opM}{M}

\DeclareMathOperator{\Aut}{Aut}

\DeclareMathOperator{\Ent}{Ent}

\DeclareMathOperator{\MA}{MA}

\DeclareMathOperator{\PSH}{PSH}

\DeclareMathOperator{\Ric}{Ric}

\DeclareMathOperator{\Spec}{Spec}

\DeclareMathOperator{\lct}{lct}

\DeclareMathOperator{\ord}{ord}

\DeclareMathOperator{\supp}{supp}

\DeclareMathOperator{\vol}{vol}

\newcommand{\NA}{\mathrm{NA}}

\renewcommand{\div}{\mathrm{div}}

\newcommand{\DNA}{\operatorname{D}}
\newcommand{\ENA}{\operatorname{E}}
\newcommand{\FNA}{\operatorname{F}}

\newcommand{\JNA}{\operatorname{J}}
\newcommand{\LNA}{\operatorname{L}}

\newcommand{\ENAstar}{\operatorname{E}^{*}}

\newcommand{\XNA}{X^{\NA}}

\newcommand{\PSHNA}{\PSH^{\NA}}
\newcommand{\cENA}{\cE^{1,\NA}}
\newcommand{\cHNA}{\cH^{\NA}}

\newcommand{\cMoneNA}{\cM^{1,\NA}}

\newcommand{\D}{\Delta}

\newcommand{\Ltwoloc}{L^2_{\mathrm{loc}}}

\numberwithin{equation}{section}       

\newtheorem{prop} {Proposition} [section]
\newtheorem{thm}[prop] {Theorem} 
\newtheorem{defi}[prop] {Definition}
\newtheorem{lem}[prop] {Lemma}
\newtheorem{cor}[prop]{Corollary}

\newtheorem{prop-def}[prop]{Proposition-Definition}

\newtheorem*{thmA}{Theorem A} 
 
\newtheorem*{thmB}{Theorem B} 
 
\newtheorem*{thmC}{Theorem C} 
\newtheorem*{thmD}{Theorem D}

\newtheorem{exam}[prop]{Example}
\newtheorem{rmk}[prop]{Remark}
\theoremstyle{remark}

\newtheorem*{ackn}{Acknowledgment} 

\begin{document}

\title{A variational approach to the Yau--Tian--Donaldson conjecture}
\date{\today}

\author{Robert J. Berman
  \and 
  S{\'e}bastien Boucksom
  \and 
  Mattias Jonsson}

\address{Mathematical Sciences\\
  Chalmers University of Technology
  and University of Gothenburg\\
  SE-412 96 G\"oteborg\\
  Sweden}
\email{robertb@chalmers.se}
\address{CNRS-CMLS\\
  \'Ecole Polytechnique\\
  F-91128 Palaiseau Cedex\\
  France}
\email{sebastien.boucksom@polytechnique.edu}
\address{Dept of Mathematics\\
  University of Michigan\\
  Ann Arbor, MI 48109--1043\\
  USA}
\email{mattiasj@umich.edu}


\begin{abstract}
 We give a variational proof of a version of the Yau--Tian--Donaldson conjecture for twisted K\"ahler--Einstein currents, and use this to express the greatest (twisted) Ricci lower bound in terms of a purely algebro-geometric stability threshold. Our approach does not involve a continuity method or the Cheeger--Colding--Tian theory, and uses instead pluripotential theory and valuations. Along the way, we study the relationship between geodesic rays and non-Archimedean metrics. 
 
\end{abstract}

\maketitle

\setcounter{tocdepth}{1}
\tableofcontents
%
%
%
%
\section*{Introduction}
The Yau--Tian--Donaldson conjecture is a central conjecture in K\"ahler geometry, whose broad goal is to relate the existence of a canonical metric in a given K\"ahler cohomology class to an algebro-geometric condition of stability. For metrics in the anticanonical class of a Fano manifold $X$, the conjecture asserts that $X$ admits a K\"ahler--Einstein metric iff $X$ is K-polystable; it was settled a few years ago by Chen--Donaldson--Sun~\cite{CDS} (see also~\cite{Tian15}), following a strategy based on a continuity method  with respect to the cone angle of a K\"ahler--Einstein metric with cone singularities along a fixed anticanonical divisor, as well as an in-depth use of the Cheeger--Colding--Tian theory of Gromov--Hausdorff limits of K\"ahler manifolds with Ricci bounds. Shortly thereafter, a proof based on the `classical' continuity method was provided by Datar and Sz\'ekelyhidi~\cite{Sze16,DaSz}, followed by another one by Chen--Sun--Wang~\cite{CSW}, based on the K\"ahler--Ricci flow. 

In the preprint version~\cite{BBJ} of the present paper, we proved that a Fano manifold $X$ without nontrivial holomorphic vector fields admits a K\"ahler--Einstein metric iff $X$ is uniformly K-stable. While only a special case of the previous results, one virtue of our approach, which is based on variational arguments and regularization techniques from pluripotential theory, lies in its relative simplicity. As we shall see, our variational method easily extends to the setting of twisted K\"ahler--Einstein currents, and contains, in particular, the smooth setting in~\cite{DaSz}, as well as the log Fano case, at least as long as the underlying variety is smooth (compare~\cite{LTW1}). Moreover, our approach naturally leads to an algebro-geometric description of the greatest (twisted) Ricci lower bound. The present paper is thus an expanded version of~\cite{BBJ}, upgraded to the setting of twisted K\"ahler--Einstein currents.

\subsection*{Main results}

Let $X$ be a (smooth) projective manifold, $L$ an ample $\Q$-line bundle on $X$, and $\theta$ a closed, quasi-positive $(1,1)$-current on $X$, \ie the sum of a positive current and a smooth form. A \emph{$\theta$-twisted K\"ahler--Einstein current} is a positive $(1,1)$-current $\om\in c_1(L)$, of finite energy in the sense of~\cite{BBGZ}, such that
\begin{equation} 
\Ric(\om)=\la\om+\theta, \,\,\,\,\,\la\in\R. \tag{TKE}
\end{equation}
When $\theta$ is a smooth form, this equation amounts to a complex Monge--Amp\`ere equation for a potential of $\om$, and pluripotential theory thus provides an interpretation of~(TKE) in the singular case as well. The constant $\la$ is determined by the cohomological condition 
$$
c_1(X,\theta):=c_1(X)-[\theta]=\la c_1(L).
$$
In the case $\la\ge 0$, the Monge--Amp\`ere formulation shows that~(TKE) only admits a solution when $\theta$ is \emph{klt}\footnote{A shorthand for Kawamata log terminal, borrowed from birational geometry.} (see Section~\ref{sec:twisted} below for details). By~\cite{BBEGZ}, each solution $\om$ of~(TKE) is an honest K\"ahler form on any open set on which $\theta$ is smooth. When $\theta$ is the integration current on an effective $\Q$-divisor $\D$ with $(X,\D)$ klt, $\om$ further has cone singularities along the normal crossing part of $\D$~\cite{GP}. 

When $\la<0$ (resp.\  $\la=0$ and $\theta$ klt), the variational approach of~\cite{BBGZ} provides a unique solution to~(TKE) (compare~\cite{BG}), generalizing classical results of Aubin and Yau~\cite{Aub,Yau}. We refer to~\cite{Tsu,ST12,Tos} for natural examples attached to Calabi--Yau fibrations, with $\theta$ of Weil-Petersson type. 

The main result of the present paper deals instead with the `twisted Fano case'. 

\begin{thmA} Let $X$ be a smooth projective manifold, $L$ an ample $\Q$-line bundle, and $\theta$ a semipositive klt current such that $c_1(X,\theta)=c_1(L)$. 
\begin{itemize}
\item[(i)] If $c_1(L)$ contains a $\theta$-twisted K\"ahler--Einstein current (resp.~a unique $\theta$-twisted K\"ahler--Einstein current), then $(X,L)$ is Ding-semistable (resp.~uniformly Ding-stable) with respect to $\theta$. 
\item[(ii)] Conversely, if $(X,L)$ is uniformly Ding-stable with respect to $\theta$, then $c_1(L)$ contains a $\theta$-twisted K\"ahler--Einstein current. 
\end{itemize}
\end{thmA}
Under a mild technical condition on the singularities of $\theta$, the twisted K\"ahler--Einstein current in (ii) is in fact unique, yielding a more symmetric statement, see Corollary~\ref{cor:coerstab}.

The notion of Ding-stability used here is phrased in terms of the non-Archimedean Ding functional, defined on test configurations and involving the log discrepancy $A_\theta(v)$ of divisorial valuations $v$ on $X$ with respect to the `klt pair' $(X,\theta)$. 
In the preprint version~\cite{BBJ} of the present paper, uniform Ding-stability was shown to be equivalent to uniform K-stability (as defined in~\cite{BHJ1,Der16}) in the usual Fano case ($L=-K_X$, $\theta=0$), building on the Minimal Model Program very much in the same way as~\cite{LX}; this equivalence was then extended to the log Fano case in~\cite{Fuj19}. As a result, we obtain: 

\begin{thmB} Let $(X,\D)$ be a \emph{log Fano manifold}, \ie $X$ is a smooth projective variety and $\D$ an effective $\Q$-divisor such that $(X,\D)$ is klt and $-(K_X+\D)$ is ample. Then the following conditions are equivalent:
\begin{itemize}
\item[(i)] $c_1(X,\D)$ contains a unique $\D$-twisted K\"ahler--Einstein current; 
\item[(ii)] $c_1(X,\D)$ contains a $\D$-twisted K\"ahler--Einstein current, and $\Aut(X,\D)$ is finite; 
\item[(iii)] the log Fano pair $(X,\D)$ is uniformly (log) K-stable. 
\end{itemize}
\end{thmB}
As already mentioned, when $\D=0$, Theorem B is basically a special case of~\cite{CDS,DaSz,CSW}. Closely related results were obtained in~\cite{LS,SW} when $\supp\D$ is a smooth divisor, and in~\cite{LTW1} in the general case (building on the preprint version~\cite{BBJ} of the present paper). 

\smallskip

As we now explain, Theorem A also yields a purely algebro-geometric description of the greatest (twisted) Ricci lower bound in terms of a stability threshold.
Given a polarized manifold $(X,L)$ and a klt current $\theta$, we define the \emph{greatest twisted Ricci lower bound} as
\begin{equation*}
  \b_\theta(X,L)=\sup\{\b\in\R\mid\Ric(\om)\ge\b\om+\theta\ 
  \text{for some $\om\in c_1(L)$}\},
\end{equation*}
where $\omega$ is a current of finite energy and where the inequality means that $\Ric(\om)-\b\om-\theta$ is smooth and semipositive.
The invariant $\b_\theta(X,L)$ is clearly bounded above by the \emph{nef threshold} $s_\theta(X,L)$, \ie the supremum of $s\in\R$ with $c_1(X,\theta)-sc_1(L)$ nef. 

In the usual smooth Fano case ($L=-K_X$, $\theta=0$), the greatest Ricci lower bound was implicitly considered in~\cite{Tian92}, and explicitly defined and studied by Rubinstein in~\cite{Rub08,Rub09}. Sz\'ekelyhidi later showed in~\cite{Sze11} that it coincides with the existence time in Aubin's continuity method. Note also that the nef threshold is equal to $1$ in that case. 

Slightly extending~\cite{BlJ,nakstab}, we introduce, on the other hand, the \emph{stability threshold} 
$$
\d_\theta(X,L):=\inf_v\frac{A_\theta(v)}{S_L(v)}, 
$$
where $v$ ranges over all divisorial valuations on $X$, $A_\theta(v)$ is the log discrepancy of $v$ mentioned above, and $S_L(v)$ is the \emph{expected vanishing order} of multisections of $L$ along $v$, defined as the limit as $m\to\infty$ of the (scaled) mean value of $v$ on sections of $mL$. Following~\cite{BlJ,nakstab}, we show that the stability threshold above coincides with the twisted analogue of the invariant originally defined in~\cite{FO}, \ie 
$$
\d_\theta(X,L)=\lim_{m\to\infty}\inf\left\{\lct_\theta(D)\mid D\text{ of $m$-basis type}\right\},
$$
where a divisor of \emph{$m$-basis type} is a $\Q$-divisor of the form 
$$
D=\frac{1}{mN_m}\sum_{j=1}^{N_m}\div(s_j)
$$ 
for some basis $(s_1,\dots,s_{N_m})$ of $H^0(mL)$, and
$$
\lct_\theta(D)=\inf_v\frac{A_\theta(v)}{v(D)}
$$
is the log canonical threshold of $D$ with respect to $\theta$. As in~\cite{nakstab}, we use non-Archimedean pluripotential theory to show that for each $\d\in\Q_{>0}$, $(X,\d L)$ is Ding-semistable (resp.\ uniformly Ding-stable) with respect to $\theta$ iff $\d\le\d_\theta(X,L)$ (resp.\ $\d<\d_\theta(X,L)$); this characterizes $\d_\theta(X,L)$ and explains the chosen terminology.

When $X$ is smooth and $\theta$ is a semipositive klt current with
$c_1(X,\theta)=c_1(L)$, Theorem~A therefore implies that the existence of a unique $\theta$-twisted K\"ahler--Einstein current is characterized by the condition $\d_\theta(X,L)>1$, at least when $\theta$ is strictly positive or has small unbounded locus, see Corollary~\ref{cor:coerstab}.

Similarly, we infer the following result from Theorem~A. 
\begin{thmC} If $(X,L)$ is a polarized manifold and $\theta$ a semipositive klt current, then 
$$
\b_\theta(X,L)=\min\{\d_\theta(X,L),s_\theta(X,L)\}.
$$ 
\end{thmC}
Here we do \emph{not} assume $c_1(X,\theta)=c_1(L)$. 
In the usual smooth Fano case, \ie $L=-K_X$ and $\theta=0$, Theorem~C was independently obtained in the appendix of~\cite{CRZ}, as a consequence of~\cite{LS,SW}, and hence ultimately~\cite{CDS} (see also~\cite{Li11} for the toric case and~\cite{Cab19} for the case of Fano manifolds of complexity one with respect to a torus action). 

The usefulness of the stability threshold to the study of K-stability has recently been further explored in several works, such as~\cite{PW,CZ,CP,BlL,BlX}. 
%
%
%
%
\subsection*{Coercivity and Ding-stability}
We now describe our strategy of proof of Theorem A. Choose a K\"ahler form $\om_0\in c_1(L)$, so that finite energy currents in $c_1(L)$ get parametrized by the space $\cE^1=\cE^1(X,\om_0)$ of finite energy potentials~\cite{GZ,BBGZ}, a complete geodesic space with respect to a metric $d_1$ introduced by Darvas~\cite{Dar15}.  

By~\cite{BBEGZ}, if $u\in\cE^1$, then $\om_u$ is a $\theta$-twisted K\"ahler--Einstein current iff $u$ minimizes the \emph{Ding functional} $\opD_\theta=\opL_\theta-\opE$, where 
$$
\opL_\theta(u)=-\tfrac12\log\int_X e^{-2u}\mu_\theta
$$
for a certain probability measure $\mu_\theta$, and $E$ is the \emph{Monge--Amp\`ere energy} functional. Further, $\opD_\theta$ admits a minimizer in $\cE^1$ as soon as it is \emph{coercive}, \ie $\opD_\theta\ge\e J-C$ for some constants $\e,C>0$, with $J\ge 0$ denoting the Aubin energy functional. A key ingredient here is the convexity of $\opD_\theta$ along plurisubharmonic (psh) geodesics in $\cE^1$, a consequence of~\cite{Bern09}. 

\smallskip

In a first step towards Theorem~A, we prove that $\opD_\theta$ is coercive iff $\opD_\theta(U_t)\to+\infty$ along each non-trivial (psh) geodesic ray $U\colon\R_{\ge 0}\to\cE^1$ (cf.~Corollary~\ref{cor:coer}), which equivalently means that the slope at infinity of $\opD_\theta$ along $U$ is positive, by convexity. The proof is based on the thermodynamical formalism of~\cite{thermo}, which shows that the coercivity of $\opD_\theta$ is equivalent to that of the twisted K-energy, and on an argument by contradiction inspired by~\cite{DaHe,DaR}, based on the entropy/energy compactness theorem of~\cite{BBEGZ} and convexity of the K-energy~\cite{BB17}. 

We next consider the set of (normal, ample) test configurations for $(X,L)$; as in~\cite{BHJ1,nakstab} we view this as a space $\cHNA$ of functions $\f$ on the set $\Xdiv$ of ($\Q$-valued) divisorial valuations on $X$. To each test configuration is attached a geodesic ray~\cite{PS,Berm16}, giving rise to a one-to-one correspondence between $\cHNA$ and geodesic rays with \emph{algebraic singularities} (emanating from $0$). It further follows from~\cite{BHJ2,Berm16} that to each functional $\opF$ among $\opE,\opL_\theta,\opD_\theta,\opJ$ corresponds a non-Archimedean version\footnote{For notational simplicity, we denote a functional and its non-Archimedean version by the same letter, dropping the superscript `NA’ used in~\cite{BHJ1,BHJ2}.} $\FNA\colon\cHNA\to\R$, with the property that 
$$
\lim_{t\to\infty} t^{-1}\opF(U_t)=\FNA(\f)
$$ 
for the geodesic ray $U$ with algebraic singularities associated to $\f\in\cHNA$. In particular, 
$$
\LNA_\theta(\f)=\inf_{\Xdiv}\left(A_\theta+\f\right),
$$
where $A_\theta>0$ denotes as above the $\theta$-twisted log discrepancy function. 

\smallskip

We say that $(X,L)$ is Ding-semistable (resp.\ uniformly Ding-stable) with respect to $\theta$ if $\DNA_\theta\ge 0$ on $\cHNA$ (resp.\  $\DNA_\theta\ge\e\JNA$ for some $\e>0$). As we just saw, uniform Ding-stability precisely means that the (Archimedean) Ding functional $\opD_\theta$ grows uniformly at infinity along geodesic rays with algebraic singularities. In order to show Theorem A, it remains to show that this condition implies that $\opD_\theta$ grows along \emph{all} non-trivial geodesic rays in $\cE^1$. 

To do this, we attach to each such ray $U$ a function $\f=U_\NA$ on $\Xdiv$, defined in terms of Lelong numbers, and compatible with the previous discussion when $U$ has algebraic singularities. Using the characterization of integrability exponents of psh functions in terms of Lelong numbers~\cite{valmul,hiro} (see Appendix~\ref{sec:valcrit}) we prove that the $\opL_\theta$ part of the Ding functional $\opD_\theta=\opL_\theta-\opE$ satisfies 
$$
\lim_{t\to\infty} t^{-1} \opL_\theta(U_t)=\LNA_\theta(U_\NA),
$$
where the right-hand side is defined by the same formula as above. On the other hand, we show that Demailly's approximation technique, based on multiplier ideals, gives rise to a sequence of rays $(U^j)$ with algebraic singularities such that $\LNA(U^j_\NA)\to\LNA(U_\NA)$ and $\ENA(U^j_\NA)\ge\lim_{t\to\infty} t^{-1}\opE(U_t)$, which is enough to conclude that $\opD_\theta(U_t)=\opL_\theta(U_t)-\opE(U_t)$ has positive slope at infinity. 

%
%
%
\subsection*{From geodesic rays to non-Archimedean functions of finite energy, and back}
As we now explain, the previous arguments admit a natural interpretation in the framework of non-Archimedean pluripotential theory, leading to a refined version of Theorem A. 

In~\cite{siminag,nama}, a non-Archimedean version of the Calabi--Yau theorem was first obtained for smooth, projective Berkovich spaces over fields of Laurent series. In~\cite{trivval}, this was adapted to the trivially valued case, in which the Berkovich analytification $\XNA$ of a projective variety $X$ provides a natural compactification of the set of divisorial valuations on $X$. Given a polarization $L$, normal, ample test configurations for $(X,L)$ are in one-to-one correspondence with non-Archimedean K\"ahler potentials, which form a space $\cHNA=\cHNA(X,L)$ of continuous functions on the compact Hausdorff space $\XNA$. Functions of finite energy are defined as decreasing limits of sequences in $\cHNA$, forming a space $\cENA$, and the non-Archimedean Calabi--Yau theorem then shows that the Monge--Amp\`ere operator induces a one-to-one correspondence between $\cENA/\R$ and Radon probability measures of finite energy on $\XNA$. 

In Section~\ref{S202}, we revisit the arguments used in the proof of Theorem A in the light of this theory. We prove that the function $U_\NA$ attached to a geodesic ray $U$ in $\cE^1$ belongs to $\cENA$, and we conversely attach to each $\f\in\cENA$ a unique maximal geodesic ray in $\cE^1$. This allows us to refine Theorem A as follows: 

\begin{thmD} Let $(X,L)$ be a polarized manifold, and $\theta$ a semipositive klt current such that $c_1(X,\theta)=c_1(L)$. The following are equivalent: 
\begin{itemize}
\item[(i)] the (Archimedean) Ding functional $\opD_\theta$ is coercive on $\cE^1$; 
\item[(ii)] the non-Archimedean Ding functional $\DNA_\theta$ is positive on all non-constant functions in $\cENA$;
\item[(iii)] $(X,L)$ is uniformly Ding-stable with respect to $\theta$. 
\end{itemize}
\end{thmD}
%
%
\subsection*{Recent developments}
Several new results have appeared since this paper was submitted, and even more since the preprint version~\cite{BBJ} first appeared. First, a version of Theorem~A was proved by Li, Tian and Wang in~\cite{LTW1,LTW2} for an arbitrary (possibly singular) log Fano pair. The proof is based on the methods in this paper, but additionally uses a delicate perturbation argument on a resolution of singularities.

As in the current paper, the existence results for K\"ahler--Einstein metrics in~\cite{LTW1,LTW2} are applicable in the case when the log Fano pair in question has finite automorphism group. The general case was recently treated by Chi Li~\cite{Li19} using a notion of uniform K-stability relative to suitable reductive subgroups of the automorphism groups. See also~\cite{His16a,His16b} for some earlier work.

Finally, there has been tremendous progress towards the Yau--Tian--Donaldson conjecture for cscK metrics. Consider a polarized manifold $(X,L)$ with finite automorphism group (modulo the scaling action of $\C^*$. The work of Chen--Cheng~\cite{CC1,CC2,CC3}, combined with~\cite{DaR,BDL}, shows that $c_1(L)$ contains a cscK metric iff the Mabuchi K-energy functional $M$ grows at infinity along each nontrivial geodesic ray in $\cE^1$. Very recently, Chi Li~\cite{Li20} proved that any geodesic ray in $\cE^1$ along which the Mabuchi energy grows at most linearly, must be a maximal geodesic ray in the sense above. This allows him to prove that if $(X,L)$ is uniformly K-stable in the sense of an inequality on $\cENA$, then $(X,L)
$ admits a cscK metric. In fact, his approach also works when the automorphism group is non-discrete.

\subsection*{Organization of the paper} 
The paper is organized as follows. 
\begin{itemize}
\item Section 1 recalls preliminary material on geodesics in the space of finite energy potentials. 
\item Section 2 reviews the thermodynamical formalism for twisted K\"ahler--Einstein currents, and proves a coercivity criterion which plays a key role in the proof of Theorem A. 
\item Section 3 discusses test configurations and Ding-stability, emphasizing the valuative point of view. 
\item Section 4 analyzes the singularities of a geodesic ray, whose Lelong numbers are encoded in a function on the space of divisorial valuations. 
\item Section 5 proves Theorem A and B above. 
\item Section 6 studies the relation between geodesic rays in $\cE^1$ and non-Archimedean functions of finite energy, and proves Theorem D. 
\item Section 7 studies the stability threshold, and proves Theorem C. 
\end{itemize}
The paper ends with two appendices where we prove certain estimates for mixed Monge--Amp\`ere integrals, and revisit the valuative criterion of integrability.

\bigskip
\begin{ackn} 
  We would like to thank Bo Berndtsson, Tam\'as Darvas, Ruadha\'\i\ Dervan, Philippe Eyssidieux, Vincent Guedj, Tomoyuki
  Hisamoto, Julius Ross, Yanir Rubinstein, Song Sun, G\'abor Sz\'ekelyhidi, Ahmed Zeriahi and David Witt Nystr\"om for helpful comments and many conversations over the years.
  We also the thank the referee for many useful remarks.
  R.B. was partially supported by the Swedish Research Council,
  the European Research Council, the Knut and Alice Wallenberg foundation, and the G\"oran Gustafsson foundation. 
  S.B. was partially supported by the ANR projects GRACK, MACK and POSITIVE\@. 
  M.J. was partially supported by NSF grants DMS-1600011 and DMS-1900025,
  the Knut and Alice Wallenberg foundation,
  and the United States---Israel Binational Science Foundation.
\end{ackn}

%
%
%
%
\section{Finite energy potentials and psh geodesics}\label{sec:E1}
In what follows, $(X,\om_0)$ denotes an $n$-dimensional compact K\"ahler manifold. In this preliminary section, we discuss plurisubharmonic (psh) paths and geodesics in the space of $\om_0$-psh functions on $X$. Most results are known, except perhaps for the characterization of geodesics given in Corollary~\ref{cor:geod}. 
%
%
\subsection{Finite energy potentials}
Denote by $\PSH:=\PSH(X,\om_0)$ the space of $\om_0$-psh functions $u:X\to[-\infty,+\infty)$, endowed with its natural weak topology, which coincides with the $L^1$-topology. The functional $u\mapsto\sup_X u$ is continuous on $\PSH$, and the space
$$
\PSH_{\sup}:=\left\{u\in\PSH\mid\sup_X u=0\right\}
$$
of sup-normalized $\om$-psh functions is compact. By~\cite{BK07}, every $u\in\PSH$ can be written as the pointwise limit of a decreasing sequence of \emph{K\"ahler potentials}, \ie elements of 
$$
\cH:=\left\{u\in C^\infty(X)\mid \om_u:=\om_0+dd^c u>0\right\}. 
$$

The \emph{Monge--Amp\`ere energy functional} $E\colon\cH\to\R$ is the antiderivative of the Monge--Amp\`ere operator $\MA(u):=V^{-1}\om_u^n$, normalized by $\opE(0)=0$. Here $V:=\int_X\om_0^n$, so that $\MA(u)$ is a probability measure. The functional $E$ is explicitly given by 
\begin{equation}\label{equ:E}
\opE(u)-\opE(v)=\frac{1}{n+1}\sum_{j=0}^n V^{-1}\int_X (u-v)\,\om_u^j\wedge\om_v^{n-j}
\end{equation}
for all $u,v\in\cH$, and hence  
\begin{equation}
  \opE(u+c)=\opE(u)+c \ \text{for $u\in\cH$, $c\in\R$};\label{e201}
\end{equation}
\begin{equation}
  u\le v\Longrightarrow \opE(u)\le \opE(v)\ \text{for $u,v\in\cH$, with equality iff $u=v$}.\label{e202}
\end{equation}
It follows that the functional $E$ admits a unique extension as a monotone, upper semicontinuous (usc) functional 
$$
\opE\colon\PSH\to\R\cup\{-\infty\},
$$
obtained by setting, for each $u\in\PSH$,
$$
\opE(u):=\inf\left\{\opE(v)\mid v\in\cH,\,v\ge u\right\}. 
$$
The space of \emph{finite energy potentials}, first introduced in~\cite{GZ} building upon the pioneering work of Cegrell~\cite{Cegrell}, can be defined as
$$
\cE^1=\cE^1(X,\om_0):=\left\{u\in\PSH\mid \opE(u)>-\infty\right\}. 
$$
We also set
$$
\cE^1_{\sup}:=\cE^1\cap\PSH_{\sup}=\left\{u\in\cE^1\mid\sup_X u=0\right\}.
$$
Unless otherwise specified, we endow $\cE^1$ with the \emph{strong topology}, defined as the coarsest refinement of the weak topology in which $E\colon\cE^1\to\R$ becomes continuous~\cite{BBGZ,BBEGZ}. 

\begin{exam} If $X$ is a Riemann surface, \ie $n=1$, a function $u\in\PSH$ belongs to $\cE^1$ iff it satisfies the classical finite energy condition $\int_X du\wedge d^c u<+\infty$, which means that the gradient of $u$ is in $L^2$. In other words, $\cE^1$ is the intersection of $\PSH$ with the Sobolev space $L^2_1$, and the strong topology is the induced Sobolev norm topology. 
\end{exam}

The following criterion for strong convergence will be useful below. 

\begin{lem}\label{lem:strongcv} A sequence $(u_j)$ in $\cE^1$ converges strongly to $u\in\cE^1$ iff $\limsup_{j\to\infty}u_j\le u$ pointwise and $\opE(u_j)\to \opE(u)$. 
\end{lem}
\begin{proof}
  Strong convergence $u_j\to u$ by definition means $\opE(u_j)\to \opE(u)$ and $u_j\to u$ weakly, and the latter property is well-known to imply 
$$
\limsup_ju_j\le(\limsup_ju_j)^*=u
$$ 
pointwise, the star denoting usc regularization.  Conversely, assume $\limsup_{j\to\infty}u_j\le u$ and $\opE(u_j)\to \opE(u)$. In order to show that $u_j\to u$ strongly, it suffices to show that the non-negative quantity
\begin{equation*}
  \opJ_u(u_j):=\int_X(u_j-u)\MA(u)+\opE(u)-\opE(u_j)
\end{equation*}
tends to $0$, by~\cite[Proposition 5.6]{BBGZ}. But this follows from Fatou's lemma, which yields
$$
\limsup_j\int_Xu_j\,\MA(u)\le\int_X(\limsup_ju_j)\MA(u)\le\int_Xu\,\MA(u),
$$
since we are dealing with functions bounded above. 
\end{proof}

By~\cite{BBEGZ}, the mixed Monge--Amp\`ere integrals 
$$
\int_X u_0\,\om_{u_1}\wedge\dots\wedge\om_{u_n}
$$ 
are well-defined for $u_0,\dots,u_n\in\cE^1$, and continuous with respect to $(u_0,\dots,u_n)$ in the strong topology. In particular,~\eqref{equ:E}--\eqref{e202} are still valid for $u,v\in\cE^1$~\cite[Theorem 4.1]{BBGZ}. 
%
%
\subsection{Psh paths}
Every connected $S^1$-invariant subset of $\C^*$ with nonempty interior is of the form 
$$
\DD_I:=\left\{\tau\in\C^*\mid -\log|\tau|\in I\right\},
$$
with $I\subset\R$ an interval (not necessarily open or closed). We are mainly interested in the case when $I$ is bounded below; then $\DD_I$ is an annulus or a punctured disc.

Slightly abusively, we will in what follows identify maps $U\colon I\to\PSH$ with $S^1$-invariant functions on $X\times\DD_I$, the correspondence being given by 
$$
U_{-\log|\tau|}(x)=U(x,\tau). 
$$
\begin{defi} A \emph{psh path}\footnote{Such a map
was called a \emph{subgeodesic} in~\cite[\S2.2]{Bern15} and subsequent works.} is a map $U\colon I\to\PSH$ defined on an open interval $I\subset\R$, such that corresponding function on $X\times\DD_I$ is $p_1^*\om_0$-psh, with $p_1\colon X\times\C\to X$ the first projection.  
\end{defi}
The condition implies in particular that $t\mapsto U_t(x)$ is convex on $I$ for each fixed $x\in X$, and hence admits limits in $[-\infty,+\infty]$ as $t$ tends to $\partial I$. Psh paths satisfy the following basic properties. 

\begin{prop}\label{prop:subcont} Let $I\subset\R$ be an open interval. Every psh path $U\colon I\to\PSH$ is continuous. If a sequence of psh paths $U^j\colon I\to\PSH$ converges to a map $U\colon I\to\PSH$ locally uniformly with respect to the $L^1$-norm, then $U$ is psh as well. 
\end{prop}
\begin{proof} Let $U\colon I\to\PSH$ be a psh path. Convexity of $t\mapsto U_t(x)$ implies that $t\mapsto\int_X U_t\,\om^n$ is convex, and hence continuous on $I$. Given $t_0\in I$, this also applies to $t\mapsto\int_X\max\{U_t,U_{t_0}\}\om^n$, as the max of two psh paths is psh. Thanks to the elementary identity
$$
\int_X|U_t-U_{t_0}|\om^n=2\int_X\left(\max\{U_t,U_{t_0}\}-U_{t_0}\right)\om^n-\int_X(U_t-U_{t_0})\om^n,
$$
we conclude that $U_t\to U_{t_0}$ in $L^1$ as $t\to t_0$, which proves the first point. 

To say that a sequence $U^j$ of psh paths converges locally uniformly to a map $U\colon I\to\PSH$ means that $U_t$ is $\om$-psh for each $t$, and $\int_X|U^j_t-U_t|\om^n$ converges to $0$ as $j\to\infty$, locally uniformly with respect to $t\in I$. By Fubini, the corresponding functions on $X\times\DD_I$ satisfy $U^j\to U$ in $L^1_{\mathrm{loc}}$. In particular, $p_1^*\om+dd^c U\ge 0$ in the sense of currents, which shows that $U$ is equal a.e.\ to an $S^1$-invariant $p_1^*\om_0$-psh function $\tilde U$ on $X\times\DD_I$. For a.e.\ $t\in I$, we thus have $U_t=\tilde U_t$ a.e.\ on $X$, and hence $U_t=\tilde U_t$ on $X$ since both functions are $\om$-psh. By local uniform convergence, the map $U\colon I\to\PSH$ is continuous. Since $\tilde U\colon I\to\PSH$ is continuous as well, and these two maps coincide outside a set of mesure $0$ in $I$, they are equal, which proves the second point. 
\end{proof}

As the next result shows, psh paths interact nicely with $\cE^1$. 
\begin{prop}\label{prop:subgeodE1} The image of any psh path $U\colon I\to\PSH$,
  with $I\subset\R$ open, is either disjoint from $\cE^1$, or contained in it. In the latter case, $U\colon I\to\cE^1$ is continuous (in the strong topology), and $t\mapsto \opE(U_t)$ is convex. 
\end{prop}

The proof relies on the following well-known computation (cf.~\cite[Proposition 6.2]{BBGZ},~\cite[\S2.4]{Bern15}).
\begin{lem}\label{lem:hessE} Assume $I\subset\R$ is open. For each smooth function $U$ on $X\times\DD_I$, the Laplacian of $\opE(U(\cdot,\tau))$ is expressed as the fiber integral
$$
dd^c_\tau \opE(U(\cdot,\tau))=V^{-1}\int_X(p_1^*\om_0+dd^c U)^{n+1}.
$$
\end{lem}

\begin{proof}[Proof of Proposition~\ref{prop:subgeodE1}] After slightly shrinking $I$, the regularization result of~\cite{BK07} yields a sequence of smooth psh paths $U^j:I\to\PSH$ decreasing pointwise to $U$. For each $j$, Lemma~\ref{lem:hessE} shows that $\opE(U^j_t)$ is a convex function of $t$. This is thus also the case for $\opE(U_t)$, which is the pointwise limit of $\opE(U^j_t)$, by continuity of $E$ along monotone sequences. By convexity of $\opE(U_t)$, the set of $t\in I$ with $\opE(U_t)=-\infty$ is either empty or equal to $I$. In the former case, we have $U_t\notin\cE^1$ for all $t$. In the latter case, the map $t\mapsto \opE(U_t)$, being convex and finite valued, is continuous on $I$, and $U\colon I\to\cE^1$ is thus continuous in the strong topology. 
\end{proof}
%
%
%
\subsection{Psh geodesics}
Following the envelope description of geodesics provided in~\cite[\S 2.2]{Bern15}, we say that a psh path $V\colon(0,1)\to\PSH$ is \emph{dominated by} two $\om$-psh functions $U_0,U_1\in\PSH$ if 
$$
\lim_{t\to 0} V_t\le U_0,\,\,\,\,\lim_{t\to 1} V_t\le U_1,
$$ 
where the pointwise limits in question exist, by convexity. If such a psh path $V$ exists, a simple envelope argument shows that there exists a largest one $U\colon (0,1)\to\PSH$, which we call the \emph{psh geodesic} joining $U_0$ to $U_1$. 

When $U_0,U_1\in\cH$ are K\"ahler potentials, X.X.~Chen's fundamental work~\cite{Che00a}, further refined in~\cite{Blo09,Blo12,CTW18}, implies that the psh geodesic joining them is $C^{1,1}$ as a function on $X\times\DD_{[0,1]}$. When $U_0,U_1$ belong to $\cE^1$, it was proved by Darvas in~\cite{Dar15} that the psh geodesic joining them exists and yields a constant speed geodesic in the Darvas metric (see \S\ref{sec:darvasmetric}). 

We provide here a direct proof of the following result, which provides an alternative characterization of psh geodesics in $\cE^1$ to be used later (see Proposition~\ref{prop:geodseq} below). 

\begin{thm}\label{thm:geod} For any pair $U_0,U_1\in\cE^1$, the psh geodesic joining them exists, and defines a continuous map $U\colon [0,1]\to\cE^1$ (in the strong topology) with $\opE(U_t)$ affine on $[0,1]$. 

Conversely, any continuous path $\tilde U\colon [0,1]\to\cE^1$ joining $U_0$ to $U_1$ with $\opE(\tilde U_t)$ affine and $\tilde U$ psh on $(0,1)$ satisfies $\tilde U=U$. 
\end{thm}

\begin{proof} Assume first that $U_0,U_1$ are bounded. As in~\cite[\S 2.2]{Bern15}, we note that for $C\gg 1$, the bounded psh path $V\colon(0,1)\to\PSH$ defined by 
$$
V_t=\max\left\{U_0-C t,U_1-C(1-t)\right\}
$$ 
is dominated by $U_0,U_1$, so the psh geodesic $U$ joining $U_0,U_1$ exists and
satisfies $V_t\le U_t$. By maximality, we have $(p_1^*\om+dd^c U)^{n+1}=0$ on $X\times\DD^*$ in the sense of pluripotential theory, and $\opE(U_t)$ is thus affine on $(0,1)$ by Lemma~\ref{lem:hessE} and a regularization argument. Further, the inequality $V_t\le U_t$ implies $\lim_{t\to 0} U_t=U_0$ and $\lim_{t\to 1} U_t=U_1$ uniformly on $X$; hence $U\colon [0,1]\to\cE^1$ is (strongly) continuous.

Let now $U_0,U_1\in\cE^1$ be arbitrary. For each $j$, denote by $U^j$ the psh geodesic joining the bounded $\om$-psh functions $U_0^j:=\max\{U_0,-j\}$ to $U_1^j:=\max\{U_1,-j\}$. Since the sequences $(U_0^j)$ and $(U_1^j)$ are decreasing, the corresponding sequence of functions $U^j$ on $X\times\DD_{[0,1]}$ is decreasing as well, thanks to the envelope description, and its limit is thus a usc function $U\colon X\times\DD_{[0,1]}\to[-\infty,+\infty)$, which is either $-\infty$ or $p_1^*\om_0$-psh on the interior $X\times\DD_{(0,1)}$. Since $\opE(U^j_t)$ is affine, we further have $\opE(U^j_t)=(1-t) \opE(U^j_0)+t \opE(U^j_1)$. By monotone continuity of $E$, it follows that $U$ induces a psh path $U\colon(0,1)\to\cE^1$ such that $\opE(U_t)=(1-t) \opE(U_0)+t \opE(U_1)$. Being usc on $X\times\DD_{[0,1]}$, it further satisfies $\lim_{t\to 0} U_t\le U_0$ and $\lim_{t\to 1} U_t\le U_1$, and Lemma~\ref{lem:strongcv} thus shows that $U\colon[0,1]\to\cE^1$ is continuous. 

Consider finally a continuous path $\tilde U\colon [0,1]\to\cE^1$ joining $U_0$ to $U_1$ with $\opE(\tilde U_t)$ affine and $\tilde U$ psh on $(0,1)$. By Lemma~\ref{lem:strongcv} again, the restriction of $\tilde U$ to $(0,1)$ is a psh path dominated by $U_0,U_1$, and hence $\tilde U_t\le U_t$ for all $t\in [0,1]$. Since $\opE(\tilde U_t)\le \opE(U_t)$ are both affine functions on $[0,1]$ with the same boundary values, they coincide, and we conclude that $\tilde U_t= U_t$. 
\end{proof}

As a direct consequence of Theorem~\ref{thm:geod}, we get: 
\begin{cor}\label{cor:geod} For a map $U\colon I\to\cE^1$ defined on a (not necessarily open, or bounded) interval, the following properties are equivalent:
\begin{itemize}
\item[(i)] the restriction of $U$ to each compact interval $[a,b]\subset I$ coincides (up to affine reparametrization) with the psh geodesic joining $U_a$ to $U_b$; 
\item[(ii)] $U$ is strongly continuous on $I$, psh on the interior $\mathring{I}$, and $\opE(U_t)$ is affine on $I$. 
\end{itemize}
\end{cor}

\begin{defi}\label{defi:geod} A map $U\colon I\to\cE^1$ satisfying the equivalent conditions of Corollary~\ref{cor:geod} is called a \emph{psh geodesic} in $\cE^1$. A \emph{psh geodesic ray} is a psh geodesic $U\colon\R_{\ge 0}\to\cE^1$. 
\end{defi}

For later use, we finally record the following mild generalization of~\cite[Theorem 1]{Dar17b} (which deals with bounded potentials). 

\begin{prop}\label{prop:supaff} Let $U\colon [a,b]\to\cE^1$ be a psh geodesic with $U_b$ more singular then $U_a$, \ie $U_b\le U_a+C$ for some constant $C>0$. Then
$$
t\mapsto\sup_X\left(U_t-U_a\right)
$$ 
is affine on $[a,b]$. In particular, if $U_a=0$ and $U_t$ is sup-normalized $($\ie $\sup_X U_t=0)$ for some $t>a$, then $U_t$ is sup-normalized for all $t\in[a,b]$.
\end{prop}
\begin{proof} After reparametrizing, we assume for ease of notation that $a=0$ and $b=1$. Set $m:=\sup_X(U_1-U_0)$. For $t\in[0,1]$, the inequality $\sup_X(U_t-U_0)\le tm$ follows directly from the convexity of $t\mapsto U_t(x)$. On the other, the psh path $V\colon(0,1)\to\PSH$ defined by $V_t=U_1+(t-1)m$ is dominated by $U_0,U_1$. By the envelope description of $U$, it follows that $U_1+(t-1)m\le U_t$ for $t\in[0,1]$, and hence 
$$
tm=\sup_X(U_1-U_0)+(t-1)m\le\sup_X(U_t-U_0),
$$
which completes the proof.
\end{proof}
%
%
\subsection{The Darvas metric}\label{sec:darvasmetric}
The weak topology of $\cE^1$ coincides with the topology induced by the $L^1(\om^n)$-norm. The strong topology of $\cE^1$, being the coarsest refinement with respect to which $E$ becomes continuous, is thus metrizable, defined by the metric
$$
d(u,v)=\|u-v\|_{L^1(\om)}+|\opE(u)-\opE(v)|.
$$
Thanks to the work of Darvas, $\cE^1$ can be equipped with a much better behaved metric. Indeed, answering a conjecture due to Guedj, it is proved in~\cite{Dar15} that $\cE^1$ can be viewed as the metric completion of $\cH$ with respect to a natural $L^1$-Finsler metric $d_1$, defined by letting $d_1(u,u')$ be the infimum of the $L^1$-lengths
$\int_0^1\|\dot u_t\|_{L^1(\MA(u_t))}dt$ 
of all smooth paths $(u_t)_{t\in[0,1]}$ in $\cH$ joining $u$ to $u'$. 

By~\cite[Corollary 4.14]{Dar15}, if $u,v\in\cE^1$ satisfy $u\ge v$, then 
$$
d_1(u,v)=\opE(u)-\opE(v). 
$$
In particular, $d_1(u,0)=-\opE(u)$ when $u\in\cE^1$ is sup-normalized. 

Finally,~\cite[Theorem 2]{Dar17a} implies that any psh geodesic $U\colon I\to\cE^1$ in the sense of Definition~\ref{defi:geod} is a constant speed geodesic for $d_1$, \ie there exists $c\ge 0$ such that 
\begin{equation}\label{e305}
  d_1(U_t,U_s)=c|t-s|
\end{equation}
for all $t,s\in I$. Note, however, that not all metric geodesics in $(\cE^1,d_1)$ are of this form.

\begin{prop}\label{prop:geodseq} If a sequence $U^j\colon I\to\cE^1$ of psh geodesics converges pointwise to a map $U\colon I\to\cE^1$, then $U$ is a psh geodesic as well. 
\end{prop}
\begin{proof} For each compact interval $[a,b]\subset I$, the $d_1$-geodesic property yields 
$$
d_1(U^j_t,U^j_s)=\left(\frac{d_1(U^j_a,U^j_b)}{|b-a|}\right)|t-s|
$$
for $t,s\in[a,b]$. It follows that $U^j$ is equicontinuous on $[a,b]$, and hence converges uniformly to $U$ on $[a,b]$, by Ascoli. As a result, $U$ is continuous on $I$, and psh on $\mathring{I}$, by Proposition~\ref{prop:subcont}. Since $U^j_t\to U_t$ strongly, the affine functions $\opE(U^j_t)$ converge pointwise to $\opE(U_t)$, which is thus affine as well, and Corollary~\ref{cor:geod} shows that $U$ is a psh geodesic. 
\end{proof}
%
%
%
%
\section{Twisted K\"ahler--Einstein currents and coercivity}\label{sec:twisted}
In this section, we review the thermodynamical formalism for twisted K\"ahler--Einstein currents, following~\cite{thermo}, and provide a coercivity criterion for certain functionals on $\cE^1$. 
%
%
\subsection{Twisted K\"ahler--Einstein currents}\label{sec:twistedKE}
In what follows, $(X,\om_0)$ denotes as above a compact K\"ahler manifold. As is well-known, smooth positive volume forms $\mu$ on $X$ are in one-to-one correspondence with Hermitian metrics $h$ on the canonical bundle $K_X$, the relation being
\begin{equation}\label{equ:mudens}
\mu=e^{-2f}i^{n^2}\Omega\wedge\bar\Omega
\end{equation}
with $f:=\log|\Omega|_h$, for any local holomorphic volume form $\Omega$. The \emph{Ricci curvature} of $\mu$ is defined as minus the curvature of $h$, \ie $\Ric(\mu)=dd^c f$ in terms of~\eqref{equ:mudens}, so that $\Ric(\om^n)=\Ric(\om)$ is the usual Ricci curvature for a K\"ahler form $\om$. 

We shall say more generally that a positive measure $\mu$ on $X$ has \emph{well-defined Ricci curvature} if it corresponds to a singular metric on $K_X$, and define its Ricci curvature $\Ric(\mu)$ as minus the corresponding curvature current. In other words, it is required that the measure $\mu$ locally satisfies~\eqref{equ:mudens} with $f\in L^1_{\loc}$, and its Ricci curvature is locally given by $\Ric(\mu)=dd^c f$. Given a closed $(1,1)$-current $\theta$, we further introduce the \emph{$\theta$-twisted Ricci curvature} of $\mu$ as 
$$
\Ric_\theta(\mu):=\Ric(\mu)-\theta,
$$
which is thus a closed $(1,1)$-current in the cohomology class
$$
c_1(X,\theta):=c_1(X)-[\theta].
$$
Note that $\Ric_\theta(\mu)$ determines $\mu$ up to a multiplicative constant. 

\begin{defi} A \emph{$\theta$-twisted K\"ahler--Einstein current} in $[\om_0]$ is a positive current of finite energy $\om=\om_u$, $u\in\cE^1$, such that $\om^n$ has well-defined Ricci curvature, and which satisfies 
\begin{equation}\label{equ:twistedKE}
\Ric_\theta(\om)=\la\om,\,\,\,\la\in\R. 
\end{equation}
\end{defi}
Here the left-hand side of~\eqref{equ:twistedKE} is defined as the twisted Ricci curvature of $\om^n$; hence $c_1(X,\theta)=\la[\om_0]$. 

\begin{lem}\label{lem:twistedKE} Assume $c_1(X,\theta)=\la[\om_0]$, let $\theta_0$ be a smooth form in the class of $\theta$, and pick a distribution $\p$ and smooth function $\rho_0$ such that $\theta=\theta_0+dd^c\p$ and $\Ric(\om_0)-\theta_0=\la\om_0+dd^c\rho_0$. For each $u\in\cE^1$, $\om=\om_u$ then satisfies~\eqref{equ:twistedKE} iff $\p\in L^1$ and 
\begin{equation}\label{equ:twistedMA}
\MA(u)=e^{2(\rho_0-\la u-\p+c)}\om_0^n
\end{equation}
for some $c\in\R$. 
\end{lem}
\begin{proof} If $\p$ is $L^1$ and $u$ solves~\eqref{equ:twistedMA}, then $\om_u$ has well-defined Ricci curvature, and 
$$
\Ric_\theta(\om_u)=\Ric\left(e^{2(\rho_0-\la u-\p+c)}\om_0^n\right)-\theta
$$
$$
=-dd^c\rho_0+\la dd^c u+dd^c\p+\Ric(\om_0)-\theta_0-dd^c\p=\la\om_u.
$$
Assume, conversely, that $\om_u$ solves~\eqref{equ:twistedKE}. Then $\om_u^n=e^{-2f}\om_0^n$ with $f\in L^1$ such that 
$$
\Ric(\om_0)+dd^c f-\theta_0-dd^c\p=\la\om_0+\la dd^c u,
$$
which implies that $f+\rho_0-\la u-\p$ is pluriharmonic on $X$, and hence constant.
\end{proof}

\begin{defi}\label{D603} We shall say that a closed $(1,1)$-current $\theta$ is \emph{klt} if $\theta$ is quasi-positive, \ie $\theta=\theta_0+dd^c\p$ with $\theta_0$ smooth and $\p$ quasi-psh, and has trivial multiplier ideal sheaf, \ie $e^{-2\p}\in L^1$. 
\end{defi}
By the solution of the openness conjecture~\cite{Bernop,GuZh}, we actually have $e^{-2\p}\in L^p$ for some $p>1$; see Appendix~\ref{sec:valcrit}.

\begin{lem}\label{lem:twistedKEklt} Let $\theta$ be a quasi-positive current, assume $c_1(X,\theta)=\la[\om_0]$ with $\la\in\R$, and let $\om\in[\om_0]$ be a $\theta$-twisted K\"ahler--Einstein current. 
\begin{itemize}
\item[(i)] If $\la\ge 0$, then $\theta$ is necessarily a klt current. 
\item[(ii)] If $\theta$ is klt, then $\om$ has continuous potentials, and is further a smooth K\"ahler form on any open set on which $\theta$ is smooth. 
\end{itemize}
\end{lem}
\begin{proof} In the notation of Lemma~\ref{lem:twistedKE}, $u$ is bounded above, and (i) thus follows directly from~\eqref{equ:twistedMA}. Assume now that $\theta$ is klt, \ie $e^{-2\p}\in L^p$ for some $p>1$. Since $u$ has zero Lelong number at each point of $X$,  a well-known result of Skoda implies that $e^{-u}$ belongs to $L^q$ for all $q<\infty$, and hence $e^{-2(\la u+\p)}\in L^{p'}$ for some $p'>1$, by H\"older's inequality. Continuity of $u$ is now a consequence of~\cite{Kolo}, while the final assertion follows from~\cite[Theorem B.1]{BBEGZ}. 
\end{proof}

\begin{exam}\label{exam:cone} If $H\subset X$ is a smooth hypersurface with iintegration current $\d_H$ and $\theta=(1-\b)\d_H$, $\b\in(0,1)$, then $\om$ is a $\theta$-twisted K\"ahler--Einstein current iff $\om$ is a (smooth) K\"ahler--Einstein metric on $X\setminus H$, with conical singularities along $H$ of cone angle $2\pi\b$. More generally, for any effective $\Q$-divisor $\D$ on $X$ with $(X,\D)$ klt, $\d_\D$-twisted K\"ahler--Einstein currents have cone singularities along the snc part of $\D$, cf.~\cite[\S 6.2]{GP}.
\end{exam}
%
%
\subsection{The Ding functional}
In what follows, we fix a klt current $\theta$, and assume that $c_1(X,\theta)=[\om_0]$. 

\begin{lem}\label{lem:muzero}There exists a unique probability measure $\mu_\theta$ such that $\Ric_\theta(\mu_\theta)=\om_0$. Further, $\mu_\theta\ge\e\om_0^n$ for some $\e>0$, and $\mu_\theta$ has $L^p$ density for some $p>1$. 
\end{lem}
\begin{proof} As noted above, any positive measure $\mu$ with well-defined Ricci curvature is uniquely determined by $\Ric_\theta(\mu)$ up to a multiplicative constant, and the uniqueness part is thus clear. To prove existence, write as above $\theta=\theta_0+dd^c\p$ and $\Ric(\om_0)-\theta_0=\om_0+dd^c\rho_0$, with $\rho_0\in C^\infty(X)$ normalized by $\int_X e^{2(\rho_0-\p)}\om_0^n=1$. Then
$$
\mu_\theta:=e^{2(\rho_0-\p)}\om_0^n
$$
yields the desired measure, which proves the final two points as well. 
\end{proof}
By Lemma~\ref{lem:twistedKE}, for each $u\in\cE^1$ we have 
$$
\Ric_\theta(\om_u)=\om_u\Longleftrightarrow\MA(u)=e^{-2u+c}\mu_\theta
$$
with $c\in\R$ a normalizing constant. 

\begin{defi}\label{defi:TL} The \emph{Ding functional} $\opD_\theta\colon\cE^1\to\R$ associated to a klt current $\theta$ such that $c_1(X,\theta)=[\om_0]$ is defined as $\opD_\theta:=\opL_\theta-\opE$, with
$$
\opL_\theta(u):=-\tfrac 12\log\int_X e^{-2u}\mu_\theta.
$$
\end{defi} 

By~\cite[\S 4]{BBEGZ} we have: 

\begin{lem}\label{lem:TDing} The Ding functional $\opD_\theta$ satisfies the following properties. 
\begin{itemize}
\item[(i)] $\opD_\theta$ is (strongly) continuous on $\cE^1$; 
\item[(ii)] if $u\in\cE^1$ minimizes $\opD_\theta$, then $\om_u$ is a $\theta$-twisted K\"ahler--Einstein current; 
\item[(iii)] if $\opD_\theta$ is coercive, then $\opD_\theta$ admits a minimizer in $\cE^1$, and $[\om_0]$ thus contains a $\theta$-twisted K\"ahler--Einstein current. 
\end{itemize}
\end{lem}
Recall that a translation invariant functional $\opF$ on $\cE^1$ is \emph{coercive} if $\opF\ge\e\opJ-C$ for some $\e,C>0$, where $\opJ\colon\cE^1\to\R_{\ge0}$ is the translation invariant functional defined by
\begin{equation}
  \opJ(u):=V^{-1}\int_Xu\,\om_0^n-\opE(u).
\end{equation}

In the semipositive case, Berndtsson's convexity results~\cite[\S 7]{Bern15} further provide: 

\begin{lem}\label{lem:postwist} If $\theta\ge 0$, then:
\begin{itemize}
\item[(i)] $\opD_\theta$ is convex along psh geodesics in $\cE^1$; 
\item[(ii)] $u\in\cE^1$ minimizes $\opD_\theta$ iff $\om_u$ is a $\theta$-twisted K\"ahler--Einstein current.
\end{itemize}
\end{lem}
%
%
%
%
\subsection{The twisted K-energy}
Consider as above a klt current $\theta$ with $c_1(X,\theta)=[\om_0]$. Note that this condition can always be achieved by choosing $\theta$ to be a smooth representative of $c_1(X)-[\om_0]$, since $\theta$ is not required to be semipositive at this stage. 

We define the \emph{$\theta$-entropy} of a probability measure $\mu$ on $X$ as (half) the entropy of $\mu$ relative to the associated probability measure $\mu_\theta$, \ie
$$
\Ent_\theta(\mu):=\tfrac 12\int_X\log\left(\frac{\mu}{\mu_\theta}\right)\mu\in[0,+\infty]
$$
if $\mu$ is absolutely continuous with respect to $\mu_\theta$, and $\Ent_\theta(\mu)=+\infty$ otherwise. It can be written as a Legendre transform 
\begin{equation}\label{equ:EntLeg}
\Ent_\theta(\mu)=\sup_{g\in C^0(X)}\left(\int g \mu-\tfrac 12\log\int e^{2g} \mu_\theta\right), 
\end{equation}
which implies that the functional $\Ent_\theta\colon\cM\to[0,+\infty]$ is convex on the space $\cM$ of probability measures, and lower continuous (lsc) in the weak topology. 

\begin{defi} The \emph{$\theta$-entropy functional} $\opH_\theta\colon\cE^1\to[0,+\infty]$ is defined by $\opH_\theta(u):=\Ent_\theta(\MA(u))$. 
\end{defi}

By~\cite[Theorem 2.17]{BBEGZ}, we have: 

\begin{lem}\label{lem:entcomp} The functional $\opH_\theta\colon\cE^1\to [0,+\infty]$ is lsc, coercive, and its sublevel sets in $\cE^1_{\sup}$ are compact in the strong topology. 
\end{lem}

\begin{defi} We say that a translation invariant functional $\opF\colon\cE^1\to\R\cup\{+\infty\}$ has \emph{$\theta$-entropy growth} if $\opF\ge \opH_\theta-A\opJ-B$ on $\cE^1$ for some constants $A,B>0$. 
\end{defi}
This condition only depends on the singularities of $\theta$. When $\theta$ is smooth, we simply say that $\opF$ has entropy growth. It then also has $\theta$-entropy growth for any klt current $\theta$, by Lemma~\ref{lem:muzero}. 

\begin{exam}\label{exam:Mab} The usual Mabuchi K-energy functional, extended to a functional 
$$
\opM\colon\cE^1\to\R\cup\{+\infty\}
$$
as in~\cite{BDL}, has entropy growth. Indeed, denoting by $\Ent(\mu)$ the entropy of a measure $\mu$ relative to $V^{-1}\om_0^n$, the Chen--Tian formula~\cite{Che00b,TianBook} expresses $\opM(u)-\Ent(\MA(u))$ as linear combination of terms of the form 
$$
\int_Xu\,\om_u^j\wedge\om_0^{n-j}\text{   and   }\int_Xu\,\Ric(\om_0)\wedge\om_u^j\wedge \om_0^{n-j-1}. 
$$
As a result, there exist $A,B>0$ with $|\opM(u)-\Ent(\MA(u))|\le A\opJ(u)+B$ (see~\eg Lemma~\ref{L302}).
\end{exam}

By Legendre duality, we have
\begin{equation}\label{equ:LEnt}
\opL_\theta(u)=\inf_{\mu\in\cM}\left(\Ent_\theta(\mu)+\int u\,\mu\right)
\end{equation}
for $u\in\cE^1$, whereas 
\begin{equation}\label{equ:EntL}
\Ent_\theta(\mu)\ge\sup_{u\in\cE^1}\left( \opL_\theta(u)-\int u\,\mu\right)
\end{equation}
for $\mu\in\cM$ (cf.~\cite[Lemma 2.11]{BBEGZ}). On the other hand, recall that the \emph{pluricomplex energy} of $\mu\in\cM$ is defined as
\begin{equation}\label{equ:EstarE}
\opE^*(\mu)=\sup_{u\in\cE^1}\left( \opE(u)-\int u\,\mu\right)\in[0,+\infty]. 
\end{equation}
Here Legendre duality shows that, for each $u\in\cE^1$,
\begin{equation}\label{equ:EEstar}
\opE(u)=\inf_{\mu\in\cM}\left(\opE^*(\mu)+\int u\,\mu\right), 
\end{equation}
the infimum being achieved precisely at $\mu=\MA(u)$. Further, the Monge--Amp\`ere operator induces a bijection between $\cE^1_{\sup}$ and the set $\cM^1\subset\cM$ of finite energy measures $\mu$~\cite{GZ,BBGZ}. 

\begin{defi} The \emph{twisted Mabuchi K-energy} $\opM_\theta\colon\cE^1\to\R\cup\{+\infty\}$ is defined by
$$
\opM_\theta(u):=\Ent_\theta(\MA(u))-\opE^*(\MA(u)).
$$ 
\end{defi}
Equivalently, 
$$
\opM_\theta(u)=\Ent_\theta(\MA(u))-\opE(u)+\int_X u\,\MA(u), 
$$
showing compatibility with the general definition of~\cite[\S 2.1]{BDL} (which does not require $\theta$ to satisfy $c_1(X,\theta)=[\om_0]$).

\begin{lem}\label{lem:legendre} The Ding and twisted Mabuchi functionals satisfy the following properties.  \begin{itemize}
     \item[(i)] $\opM_\theta$ has $\theta$-entropy growth.
     \item[(ii)] $\opM_\theta(u)\ge \opD_\theta(u)$ for all $u\in\cE^1$, with equality iff $\om_u$ is a $\theta$-twisted K\"ahler--Einstein current.
      \item[(iii)] $\inf_{\cE^1} \opM_\theta=\inf_{\cE^1}\opD_\theta\in\R\cup\{-\infty\}$; in particular, $\opD_\theta$ is bounded below iff $\opM_\theta$ is. 
     \item[(iv)] $\opD_\theta$ is coercive iff $\opM_\theta$ is. 
     \item[(v)] If $\theta\ge 0$, then $\opM_\theta$ is geodesically convex on $\cE^1$. 
  \end{itemize}
\end{lem}
\begin{proof} (i)--(iii) are proved just as in~\cite{thermo,BBEGZ}. Indeed, (i) is a consequence of the known estimates $n^{-1}\opJ(u)\le \opE^*(\MA(u))\le n\opJ(u)$; (ii) follows from~\eqref{equ:LEnt}--\eqref{equ:EEstar}, and implies $\inf_{\cE^1} \opM_\theta\ge\inf_{\cE^1}\opD_\theta$. To prove the reverse inequality, we can assume that $c:=\inf \opM_\theta>-\infty$. Then $\Ent_\theta(\mu)\ge c+\opE^*(\mu)$ for all $\mu$, and hence 
$$
\opL_\theta(u)=\inf_\mu\left(\Ent_\theta(\mu)+\int u\,\mu\right)\ge c+\inf_\mu\left(\opE^*(\mu)+\int u\,\mu\right)=c+\opE(u),
$$
\ie $\opD_\theta\ge c$, which proves (iii). Next, (v) follows from~\cite[Theorem 1.2]{BDL}, itself a consequence of~\cite{BB17,CLP}. It remains to prove~(iv), for which we argue as in~\cite[Corollary 3.6]{thermo}. Since $\opM_\theta\ge \opD_\theta$, $\opM_\theta$ is coercive as soon as $\opD_\theta$ is.
  For the converse, the surjectivity of the Monge--Amp\`ere operator
  $\MA\colon\cE^1\to\cM^1$ implies that the coercivity of $\opM_\theta$ is equivalent to 
  existence of $C>0$ and $\e\in(0,1)$ such that $\opE^*(\mu)\le\e \Ent_\theta(\mu)+C$ for
  all probability measures $\mu$.
  (Recall that $\opE^*(\mu)=\infty$ implies $\Ent_\theta(\mu)=\infty$, cf.~\cite[Lemma 2.18]{BBEGZ}).
  By~\eqref{equ:EEstar} and~\eqref{equ:LEnt}, we infer
  $$
  \opE(\e u)=\inf_\mu\left(\opE^*(\mu)+\int \e u\,\mu\right)\le\e\inf_\mu\left(\Ent_\theta(\mu)+\int u\,\mu\right)+C=\e \opL_\theta(u)+C.
  $$
  Normalizing $u$ by $\int u\,\om^n=0$, an inequality due to Ding~\cite[Remark 2]{Din88} yields
  $$
  -\opE(\e u)=\opJ(\e u)\le\e^{1+\tfrac 1n}\opJ(u)=-\e^{1+\tfrac 1n} \opE(u). 
  $$
  We infer $\e' \opE(u)\le \opL_\theta(u)+C'$ with $\e':=\e^{1/n}$ and $C'=\e^{-1}C$, which gives the coercivity estimate
  $$
  \opD_\theta(u)=\opL_\theta(u)-\opE(u)\ge(\e'-1)\opE(u)-C'=(1-\e')\opJ(u)-C'
  $$
  since $\e'<1$ and $\int u\,\om^n=0$. 
 \end{proof}
%
%
%
\subsection{The coercivity criterion}
The next result is based on the first version of the present paper~\cite[\S 2.3]{BBJ}, itself inspired by~\cite{DaHe,DaR}. The statement below is basically~\cite[Theorem 3.6]{ICM}, but see also~\cite[Theorem 6.1]{CC2} for a closely related result. 

\begin{thm}\label{thm:coer} Let $\opF\colon\cE^1\to\R\cup\{+\infty\}$ be a translation invariant functional, and assume that $\opF$ is lsc, geodesically convex, and has $\theta$-entropy growth for some klt current $\theta$ such that $c_1(X,\theta)=[\om_0]$. 
\begin{itemize}
\item[(a)] If $\opF$ is coercive, then it admits a minimizer in $\cE^1$. 
\item[(b)] If $\opF$ is not coercive, then given any $u\in\cE^1$ there exists a nontrivial psh geodesic ray $U\colon\R_{\ge 0}\to\cE^1$ emanating from $u$ along which $\opF(U_t)$ decreases. 
\end{itemize}
\end{thm}
Here we say that $U$ is trivial if $U_t-U_0$ only depends on $t$.

\begin{exam} The dichotomy (a)--(b) applies in particular to the Mabuchi K-energy functional of any compact K\"ahler manifold, which is geodesically convex by~\cite{BB17,CLP,BDL}, and has entropy growth by Example~\ref{exam:Mab}.
\end{exam}

\begin{proof}[Proof of Theorem~\ref{thm:coer}] Assume that $\opF$ is coercive, and pick a minimizing sequence $(u_j)$ in $\cE^1$, such that $\lim_{j\to\infty}\opF(u_j)=\inf \opF>-\infty$. By translation invariance, we may assume $u_j$ is sup-normalized. Since $\opF(u_j)$ is bounded above, so is $\opJ(u_j)$, by coercivity, and the entropy growth assumption implies that $\opH_\theta(u_j)$ is bounded as well. By Lemma~\ref{lem:entcomp}, $(u_j)$ stays in a (strongly) compact subset of $\cE^1$, and may thus be assumed to converge. Since $\opF$ is lsc, the limit is then a minimizer of $\opF$. 

Assume, conversely, that $\opF$ is not coercive, and pick sequences $u_j\in\cE^1_{\sup}$, $\e_j\searrow 0$ and $C_j\to+\infty$ such that 
\begin{equation}\label{equ:Fcoer}
\opF(u_j)\le\e_j \opJ(u_j)-C_j.
\end{equation}
By entropy growth, we have $\opF(u_j)\ge -A \opJ(u_j)-B$ for some constants $A,B>0$; hence $(A+\e_j) \opJ(u_j)\ge C_j-B$, which shows that 
$$
T_j:=d_1(u_j,0)=-\opE(u_j)=\opJ(u_j)+O(1)
$$ 
tends to $+\infty$. Denote by $U^j\colon[0,T_j]\to\cE^1$ the unit speed psh geodesic connecting $u$ to $u_j$; this takes values in $\cE^1_{\sup}$ by Proposition~\ref{prop:supaff}. By convexity of $\opF$ along $U^j$, we get for, $j\gg 1$ and all $t\in [0,T_j]$, 
\begin{equation}\label{equ:Mconv}
\opF(U^j_t)-\opF(u)\le t T_j^{-1}(\opF(u_j)-\opF(u))\le t\e_j. 
\end{equation}
By $\theta$-entropy growth of $\opF$, it follows that the $1$-Lipschitz maps $U^j\colon[0,T_j]\to(\cE^1,d_1)$ send every given compact subset of $\R_{\ge 0}$ to a fixed subset of $\cE^1_{\sup}$ with bounded $\theta$-entropy, and hence compact for the metric space topology, by Lemma~\ref{lem:entcomp}. By the general Arzel\`a-Ascoli theorem (for maps into metric spaces), $U^j$ therefore converges uniformly on compact sets of $\R_{\ge 0}$ to a continuous map $U\colon\R_{\ge 0}\to\cE^1_{\sup}$, after perhaps passing to a subsequence. By Proposition~\ref{prop:geodseq}, $U$ is a psh geodesic ray, and $\opF(U_t)\le\opF(u)$ by~\eqref{equ:Mconv} and lower semicontinuity; this implies that $t\mapsto\opF(U_t)$ is decreasing, by convexity. 
\end{proof}

\begin{cor}\label{cor:coer} Assume $\theta\ge 0$. If $\opM_\theta$, or equivalently $\opD_\theta$, is not coercive, then given any $u\in\cE^1_{\sup}$ there exists a nonconstant psh geodesic ray $U\colon\R_{\ge 0}\to\cE^1_{\sup}$ emanating from $u$ such that $\opD_\theta(U_t)\le \opM_\theta(U_t)\le \opM_\theta(u)$. 
\end{cor}
\begin{proof} The lsc functional $\opM_\theta\colon\cE^1\to\R\cup\{+\infty\}$ has $\theta$-entropy growth, and is geodesically convex by Lemma~\ref{lem:legendre}, since $\theta\ge 0$. The result thus follows from Theorem~\ref{thm:coer}. 
\end{proof}

As a further consequence of Theorem~\ref{thm:coer}, we obtain the following version of~\cite[Theorem 2.12]{DaR}.

\begin{thm}\label{thm:Tboundcoer} Let $\theta$ be a semipositive klt current with $c_1(X,\theta)=[\om_0]$.  
\begin{itemize}
\item[(i)] If $[\om_0]$ contains a $\theta$-twisted K\"ahler--Einstein current, then $\opD_\theta$ and $\opM_\theta$ are bounded below on $\cE^1$; 
\item[(ii)] if $[\om_0]$ contains a \emph{unique} $\theta$-twisted K\"ahler--Einstein current, then $\opD_\theta$ and $\opM_\theta$ are coercive on $\cE^1$. 
\end{itemize}
\end{thm}

\begin{proof} Assume given $u\in\cE^1$ with $\Ric_\theta(\om_u)=\om_u$. By Lemma~\ref{lem:postwist}, $\inf_{\cE^1} \opD_\theta=\opD_\theta(u)>-\infty$, and $\opM_\theta$ is bounded below as well, by Lemma~\ref{lem:legendre}. If $\opM_\theta$ fails to be coercive, Corollary~\ref{cor:coer} yields a non-constant psh geodesic ray $U\colon\R_{\ge 0}\to\cE^1_{\sup}$ emanating from $u$ such that $\opM_\theta(U_t)\le \opM_\theta(u)$ for all $t\ge 0$. Using Lemma~\ref{lem:legendre} again, we infer 
$$
\opD_\theta(U_t)\le \opM_\theta(U_t)\le \opM_\theta(u)=\opD_\theta(u)=\inf_{\cE^1} \opD_\theta.  
$$
By Lemma~\ref{lem:TDing}, $\om_{U_t}$ provides a whole ray of twisted K\"ahler--Einstein currents in $c_1(X,\theta)$.
\end{proof}
%
%
%
\subsection{Uniqueness of twisted K\"ahler--Einstein currents}\label{sec:unique}
We briefly discuss here how coercivity implies uniqueness of twisted K\"ahler--Einstein currents, under a mild regularity assumption on $\theta$. 

\begin{lem}\label{lem:uniqueKE} Let $\theta$ be a semipositive klt current such that $c_1(X,\theta)=[\om_0]$, and assume that one of the following two conditions holds:
\begin{itemize}
\item[(i)] $\theta$ has small unbounded locus, \ie its local potentials are locally bounded outside a closed complete pluripolar subset of $X$. 
\item[(ii)] $\theta$ is strictly positive, \ie $\theta\ge\e\om_0$ for some $\e>0$. 
\end{itemize}
If $\opD_\theta$ (or, equivalently, $\opM_\theta$) is coercive, then $[\om_0]$ contains a unique $\theta$-twisted K\"ahler--Einstein current. 
\end{lem}
\begin{proof} Existence is already a consequence of Lemma~\ref{lem:TDing}. In case (ii), uniqueness follows from the strict convexity of the twisted K-energy, as in~\cite[Theorem 4.13]{BDL}. In case (i),~\cite[Theorem 6.1]{Bern15} applies, and shows that the geodesic segment joining any two $\theta$-twisted K\"ahler--Einstein currents $\om,\om'$ is realized by the flow $\om_t:=\exp(tv)^*\om$ of a holomorphic vector field $v$. Now $\om_t$ makes sense for all $t\in\R$, extending the geodesic segment to a whole geodesic line, and we can then argue as in the proof of~\cite[Theorem 5.4]{BBEGZ}. Indeed, the coercivity property implies that $\opJ_\om(\om_t)$ is bounded, the Ding functional being constant along $\om_t$. By convexity, $\opJ_\om(\om_t)$ is thus constant, and hence $\om_t=\om$ for all $t$, yielding in particular $\om'=\om_1=\om$. 
\end{proof}

%
%
%
\section{Valuations and stability}\label{sec:KDingstab}
In this section, $X$ is smooth complex projective variety endowed with an ample $\Q$-line bundle $L$. We use~\cite{BHJ1,trivval,nakstab} as references. 
%
%
\subsection{Log discrepancy}\label{S302}
Denote by $\Xdiv$ the set of (rational) \emph{divisorial valuations} on $X$, \ie valuations $v\colon\C(X)^*\to\Q$ of the form $v=c\ord_E$ with $c\in\Q_{>0}$ and $E$ a prime divisor on some normal variety $Y$ mapping birationally to $X$. The \emph{log discrepancy} of $v\in\Xdiv$ is 
$$
A_X(v):=c\left(1+\ord_E(K_{Y/X})\right),
$$
where $K_{Y/X}$ denotes the relative canonical divisor. It is convenient to also include in $\Xdiv$
the \emph{trivial valuation} on $\C(X)$. This will be denoted by $v_\triv$, and has $A_X(v_\triv)=0$. 

The projection $p_1\colon X\times\C\to X$ induces a map
$(X\times\C)^{\mathrm{div}}\to\Xdiv$; this has a canonical section $\sigma\colon\Xdiv\to(X\times\C)^{\mathrm{div}}$, the \emph{Gauss extension}, defined by 
$$
\sigma(v)\left(\sum_i f_i\tau^i\right)=\min_i\left\{v(f_i)+i\right\}
$$
for each finite sequence of functions $f_0,\dots,f_r\in \C(X)$, with $\tau$ denoting the coordinate on $\C$. The image of $\sigma$ consists precisely of rational divisorial valuations $w$ on $X\times\C$ that are
$\C^*$-invariant (under the action on the second factor), and normalized
by $w(\tau)=1$. For each $v\in\Xdiv$ we have 
$$
A_{X\times\C}(\sigma(v))=A_X(v)+1.
$$ 
Assume now that we are given a quasi-positive closed $(1,1)$-current $\theta$, and write as above $\theta=\theta_0+dd^c\p$ with $\theta_0$ smooth and $\p$ quasi-psh. For each $v\in\Xdiv$ we can make sense of $v(\theta)=v(\p)$ as a generic Lelong number on some blowup, see~\cite{pshsing,hiro} and Appendix~\ref{sec:valcrit}.

\begin{defi} The \emph{$\theta$-twisted log discrepancy function} $A_\theta\colon\Xdiv\to\R$ is defined by setting 
$$
A_\theta(v):=A_X(v)-v(\theta).
$$
\end{defi}

\begin{exam} When $\theta$ is smooth, $A_\theta$ is simply equal to $A_X$. When $\theta=\d_\D$ is the integration current of a $\Q$-divisor $\D$, $A_\theta=A_{(X,\D)}$ is the usual log discrepancy of the pair $(X,\D)$. 
\end{exam} 
As explained in Appendix~\ref{sec:valcrit}, it follows from the openness conjecture and the valuative criterion of integrability that $\theta$ is klt (in the sense of Definition~\ref{D603}) iff there exists $\e>0$ such that $A_\theta\ge\e A_X$ on $\Xdiv$; see Corollary~\ref{cor:lctval1}.
%
%
\subsection{Test configurations and non-Archimedean potentials}\label{S201}
Recall that a \emph{test configuration} $(\cX,\cL)$ for $(X,L)$ is a $\C^*$-equivariant partial compactification over $\C$ of $(X,L)\times\C^*$; more precisely, it consists of a flat projective morphism $\pi\colon\cX\to\C$, a $\Q$-line bundle $\cL$ on $\cX$, a $\C^*$-action on $(\cX,\cL)$ lifting the standard one on $\C$, and an identification of the fiber over $1\in\C$ with $(X,L)$. We say that $(\cX,\cL)$ is normal (resp.~ample) when $\cX$ is normal (resp.~$\cL$ is relatively ample). 

Each test configuration $(\cX,\cL)$ defines a non-Archimedean metric on the Berkovich analytification of $L$ with respect to the trivial absolute value on $\C$; this will be viewed in the present paper through its canonical potential $\f=\f_{(\cX,\cL)}$, a function on $\Xdiv$ defined as follows. Pick a test configuration $\cX'$ dominating both $\cX$ and the trivial test configuration $X\times\C$, with $\C^*$-equivariant morphisms $\rho\colon\cX'\to\cX$ and $\mu\colon\cX'\to X\times\C$, and let $p_1\colon X\times\C\to X$ be the projection. Then $\rho^*\cL=\mu^*p_1^*L+D$ for 
a unique $\Q$-divisor $D$ supported on the central fiber, and we set, for each $v\in\Xdiv$, 
$$
\f(v):=\sigma(v)(D)
$$
with $\sigma(v)$ the $\C^*$-invariant lift of $v$ as above. 

The trivial test configuration induces the zero function. Two test configurations $(\cX,\cL)$ and $(\cX',\cL')$ determine the same function on $\Xdiv$ iff the pullbacks of $\cL$ and $\cL'$ to some test configuration dominating $\cX$ and $\cX'$ coincide. Further, the map $(\cX,\cL)\mapsto\f_{(\cX,\cL)}$ is injective on the set of normal, ample test configurations. Its image is denoted by $\cHNA$. Functions attached to arbitrary test configurations are then differences of functions in $\cHNA$. 

Functions in $\cHNA$ can alternatively be described in terms of $\C^*$-invariant ideals on $X\times\C$ (called \emph{flag ideals} in~\cite{Oda13}). Denoting by $\tau$ the coordinate on the $\C$-factor, each such ideal is of the form $\fa=\sum_{i=0}^r\tau^i\fa_i$ for a sequence of ideals $\fa_0\subset\dots\subset\fa_r$ on $X$, and defines a function $\f_\fa$ on $\Xdiv$ by setting for $v\in \Xdiv$ 
$$
\f_\fa(v):=-\sigma(v)(\fa)=\max_i\{-v(\fa_i)-i\}. 
$$
A function $\f\colon\Xdiv\to\R$ then belongs to $\cHNA$ iff it is of the form $\f=m^{-1}\f_\fa+c$ with $c\in\Q$, $m\in\N^*$ and $\fa$ a $\C^*$-invariant ideal on $X\times\C$, cosupported on $X\times\{0\}$ (\ie $\fa_r=\cO_X$ in the above notation), and such that the sheaf $p_1^*(mL)\otimes\fa$ is globally generated on $X\times\C$ (\ie $mL\otimes\fa_i$ globally generated for all $i$). Using this description, it is easy to check: 
\begin{lem}\label{lem:nabounded} Each function $\f\in\cHNA$ is bounded on $\Xdiv$, with $\sup_{\Xdiv}\f=\f(v_\triv)$. 
\end{lem}

%
%
\subsection{Non-Archimedean functionals and stability}\label{S204}
In~\cite{BHJ1}, non-Archimedean versions of a number of usual functionals on $\cH$ were introduced. They are defined as functionals on $\cHNA$, the idea being that the non-Archimedean version of a functional $F$ should compute the slopes at infinity of $F$ along psh rays in $\cHNA$ with algebraic singularities in the sense of \S\ref{sec:algsing}.\footnote{As opposed to the convention in~\cite{BBJ,BHJ1} we do not include ``NA'' in the notation for the non-Archimedean functionals. However, these functionals are defined on $\cHNA$ rather than $\cH$.} 

First, we define non-Archimedean versions of the Monge--Amp\`ere energy and $J$-energy by 
\begin{equation}\label{equ:ENA}
  \ENA(\f)=\frac{\left(\bar\cL^{n+1}\right)}{(n+1)V}
  \quad\text{and}\quad
  \JNA(\f)=\sup\f-\ENA(\f)
\end{equation}
for all $\f\in\cHNA$, where $(\bar\cX,\bar\cL)$ is the compactification of the unique normal, ample test configuration $(\cX,\cL)$ such that $\f=\f_{(\cX,\cL)}$, and $V=(L^n)$ is the volume of $L$.

\smallskip
Now fix a klt current $\theta$ on $X$. 
\begin{defi}\label{defi:NADing} The \emph{non-Archimedean Ding functional} $\DNA_\theta\colon\cHNA\to\R$ with respect to $\theta$ is defined as $\DNA_\theta:=\LNA_\theta-\ENA$ with 
\begin{equation}\label{eqn:DingNA}
  \LNA_\theta(\f)=\inf_{\Xdiv}(A_\theta+\f). 
\end{equation}
\end{defi} 
Since $A_\theta\ge 0$, $\LNA_\theta(\f)\ge\inf\f$ is indeed finite, by Lemma~\ref{lem:nabounded}. Note also that non-Archimedean Ding functional $\DNA_\theta$, in contrast to its Archimedean counterpart, only depends on the singularities of $\theta$, and thus makes sense without requiring $c_1(X,\theta)=c_1(L)$ (which could anyway always be achieved by adding a smooth form to $\theta$). 

\begin{defi}\label{defi:Dingstab} The polarized variety $(X,L)$ is \emph{Ding-semistable} (resp.\  \emph{uniformly Ding-stable}) with respect to $\theta$ if $\DNA_\theta\ge 0$ on $\cHNA$ (resp.\  $\DNA_\theta\ge\e
\JNA$ on $\cHNA$ for some $\e>0$). 
\end{defi}
When $\theta$ is smooth, $A_\theta=A_X$, and we thus drop the reference to $\theta$ in the above definitions, as in~\cite{nakstab}. By~\cite[Corollary 2.11, Theorem 2.12]{nakstab}, Ding-semistability (resp.~uniform Ding-stability) of $(X,L)$ implies (and is conjecturally equivalent to) twisted K-semistability (resp.~uniform twisted K-stability) in the twisted Fano case, in the sense of~\cite{Der16}. 

\smallskip

Suppose now that $\theta$ is the integration current on an effective $\Q$-divisor $\D$ with $(X,\D)$ klt and $c_1(L)=c_1(X,\D)$. Ding-stability with respect to $\theta$ then coincides with Ding-stability of the log Fano variety $(X,\D)$, as studied in~\cite{BHJ1,Fuj19}, and we thus have: 

\begin{thm}\label{thm:DvsK}\cite{BBJ,Fuj19} Let $\D$ be an effective $\Q$-divisor with $(X,\D)$ klt and $c_1(L)=c_1(X,\D)$. Then $(X,L)$ is Ding-semistable (resp.\  uniformly Ding-stable) with respect to $\D$ iff the log Fano variety $(X,\D)$ is log K-semistable (resp.~uniformly log K-stable). 
\end{thm}
When $\D=0$, this result was indeed proved in the preprint version~\cite{BBJ}, the  argument relying on the Minimal Model Program along the lines of~\cite{LX}. The result was later extended to the general log Fano case in~\cite{Fuj19}. 

%
%
%
%
\section{Psh rays and Lelong numbers}
In this section we study psh rays of linear growth, to which we associate functions on $\Xdiv$ defined in terms of Lelong numbers. We also introduce the class of rays with algebraic singularities; these provide a bridge between psh rays and test configurations.

In what follows, we fix an ample $\Q$-line bundle on $X$ and a K\"ahler form $\om_0\in c_1(L)$.
%
%
%
%
\subsection{Rays of linear growth}
For each psh ray $U\colon\R_{>0}\to\PSH$, $\sup_X U_t$ is a convex function of $t$. As a result, $\sup_X U_t\ge -Ct$ for some $C>0$ as $t\to\infty$, and the slope at infinity
\begin{equation}\label{equ:lmax}
\la_{\max}:=\lim_{t\to\infty} t^{-1}\sup_X U_t
\end{equation}
exists in $\R\cup\{+\infty\}$. We have $\la_{\max}<\infty$ iff $\sup_X U_t=O(t)$, in which case we say that $U$ has \emph{linear growth}. For rays in $\cE^1$, we equivalently have: 

\begin{prop}\label{prop:georaylin}
  A psh ray $U\colon\R_{>0}\to\cE^1$ has linear growth iff $d_1(U_t,0)=O(t)$ as $t\to\infty$. In particular, any psh geodesic ray has linear growth. 
\end{prop}
\begin{proof} By Proposition~\ref{prop:subgeodE1}, $\opE(U_t)$ is convex, and hence admits a linear lower bound $\opE(U_t)\ge -C t$ for $t\ge 1$. Assume $U$ has linear growth, and pick $a>0$ such that $U_t\le at$ for $t\ge 1$. Then $d_1(U_t,at)=at-\opE(U_t)\le C't$, and $d_1(U_t,0)=O(t)$, by the triangle inequality. Assume, conversely, that $d_1(U_t,0)=O(t)$. By~\cite[Proposition 2.7]{GZ1}, 
$$
\sup_X U_t=V^{-1}\int_X U_t\om^n+O(1),
$$
while Corollary~\ref{cor:intd} in the appendix gives $\left|\int_X U_t\,\om^n\right|\le C_n d_1(U_t,0)$; hence $\sup_XU_t=O(t)$ and the result follows.
\end{proof}
%
%
%
\subsection{Lelong numbers}
For a psh ray $U\colon\R_{>0}\to\PSH$ of linear growth, $U-a t$ is bounded above as $t\to\infty$, for some $a\in\R$. Equivalently, the $S^1$-invariant $p_1^*\om_0$-psh function $V$ on $X\times\DD^*$ defined by
$$
V(x,\tau):=U_{-\log|\tau|}(x)+a\log|\tau|
$$
is bounded above near $X\times\{0\}$, and hence uniquely extends to a quasi-psh function on $X\times\DD$. For each divisorial valuation $w$ on $X\times\C$, we can make sense of $w(V)\ge 0$ as a generic Lelong number on a suitable blowup, see~\cite{hiro} and Appendix~\ref{sec:valcrit}. Following~\cite[\S5]{Bermzero}, we set $w(U):=w(V)-a w(\tau)$; this is independent of the choice of $a$ by additivity of Lelong numbers. 

\begin{defi}\label{defi:Lelong} To each psh ray $U\colon\R_{>0}\to\PSH$ of linear growth  we associate a function 
$$
U_\NA\colon\Xdiv\to\R
$$ 
by setting $U_\NA(v)=-\sigma(v)(U)$ for $v\in\Xdiv$. 
\end{defi}
Recall that $\sigma\colon\Xdiv\to(X\times\C)^{\mathrm{div}}$ denotes Gauss extension, cf.~\S\ref{S302}. It sends the trivial valuation $v_\triv$ to $\ord_{X\times\{0\}}$. 
\begin{lem}\label{lem:slopesup} We have $U_\NA(v_\triv)=\sup_{\Xdiv} U_\NA=\la_{\max}$ (see~\eqref{equ:lmax}).
\end{lem}
\begin{proof} After adding a linear function of $t$, we may assume that $U$ itself extends to a quasi-psh function on $X\times\DD$. The left-hand side is then minus the generic Lelong number of $U$ along $X\times\{0\}$, which is also the maximum of all $c\ge 0$ such that $U\le c\log|\tau|+O(1)$ near $X\times\{0\}$, \ie $\sup_X U_t\le-c t+O(1)$. By convexity of $t\mapsto\sup_X U_t$, we infer $U_\NA(v_\triv)=\lim_{t\to\infty} t^{-1}\sup_X U_t$. Finally, if $U\le c\log|\tau|+O(1)$ for some $c\ge 0$, then $w(U)\ge c w(\tau)$ for every divisorial valuation $w$ on $X\times\C$, and hence $U_\NA(v)\le U_\NA(v_\triv)$ for all $v\in\Xdiv$. 
\end{proof}

%
%
%
\subsection{Relation to the Ross--Witt Nystr\"om Legendre transform}
Recall from~\cite[\S 6]{RWN} that the \emph{Legendre transform} of a psh ray $U\colon\R_{>0}\to\PSH$ is the concave family of functions $(\hat U^\la)_{\la\in\R}$ on $X$ defined by
$$
\hat U^\la:=\inf_{t>0}\left\{U_t-t\la\right\}.
$$
By the Kiselman minimum principle, for each $\la$ we either have $\hat U^\la\in\PSH$ or $\hat U^\la\equiv-\infty$, and Legendre duality yields
\begin{equation}\label{equ:legray}
U_t=\sup_{\la\in\R}\left\{\hat U^\la+\la t\right\}.
\end{equation}
By~\eqref{equ:legray}, $\sup_X U_t=\sup_\la\left\{\sup_X \hat U^\la+\la t\right\}$, which shows that
\begin{equation}\label{equ:maxleg}
\la_{\max}=\sup\left\{\la\in\R\mid \hat U^\la\ne-\infty\right\}\text{     and     }U_t=\sup_{\la<\la_{\max}}\left\{\hat U^\la+\la t\right\}. 
\end{equation}

Assuming $U$ is of linear growth, \ie $\la_{\max}<\infty$, the function $U_\NA$ can also be described in terms of the Legendre transform $(\hat U^\la)$. For each $\la<\la_{\max}$, $\hat U^\la$ is a quasi-psh function on $X$, and we can thus define $\hat U^\la_\NA\colon\Xdiv\to\R_{\le 0}$ by $\hat U^\la_\NA(v):=-v(\hat U^\la)$. This function is homogeneous of degree $1$ with respect to the scaling action of $\Q_{>0}$, and we have
$$
U_\NA=\sup_{\la<\la_{\max}}\left\{\hat U^\la_\NA+\la\right\}.
$$

%
%
%
%
\subsection{Algebraic singularities}\label{sec:algsing}
Choosing a smooth Hermitian metric $h_0$ on $L$ with curvature $\om_0$ sets up a one-to-one correspondence between psh rays $U\colon\R_{>0}\to\PSH$ and $S^1$-invariant psh metrics $e^{-2U} p_1^*h_0$ on $(X\times\DD^*,p_1^*L)$. We say that $U$ \emph{induces a psh metric} on a normal test configuration $(\cX,\cL)$ if the corresponding psh metric on $(X\times\DD^*,p_1^*L)\simeq(\cX,\cL)|_{\DD^*}$ extends to a psh metric on $(\cX,\cL)|_{\DD}$. 

\begin{lem}\label{lem:pshtest} Given a psh ray $U\colon\R_{>0}\to\PSH$ and a normal test configuration $(\cX,\cL)$, the following conditions are equivalent: 
\begin{itemize}
\item[(i)] $U$ induces a psh metric on $(\cX,\cL)$;
\item[(ii)] $U$ has linear growth, and $U_\NA\le\f_{(\cX,\cL)}$. 
\end{itemize}
If the induced psh metric in (i) is further locally bounded, then $U_\NA=\f_{(\cX,\cL)}$. 
\end{lem} 
\begin{proof} By normality of $\cX$, a psh metric on $(\cX,\cL)|_{\DD^*}$ extends to $(\cX,\cL)|_\DD$ iff the same holds for its pull-back to a higher test configuration. After passing to a higher test configuration, we may thus assume, without loss of generality,  that $\cX$ is smooth and dominates the trivial test configuration via $\rho\colon\cX\to X\times\C$. Write $\cL=\rho^*p_1^*L+D$, and pick a positive integer $m$ such that $mD$ is Cartier. Then (i) holds iff $U+m^{-1}\log|f|$ is locally bounded above for any choice of local equation $f$ for $mD$. Since $D+a\cX_0$ is effective for $a>0$ large enough, it follows that $U_t\le at+O(1)$, which shows that $U$ has linear growth. For any divisorial
  valuation $w$ on $X\times\C$ with $w(\tau)>0$, we also get $w(U)\ge -m^{-1} w(f)=-w(D)$. Applying this to $w=\sigma(v)$ with $v\in\Xdiv$ shows that $U_\NA(v)\le\sigma(v)(D)=\f_{(\cX,\cL)}(v)$. This proves (i)$\Longrightarrow$(ii), and the final assertion is proved similarly. 
  
Conversely, assume (ii). Then $\ord_E(U)\ge-\ord_E(D)$ for each irreducible component $E$ of $\cX_0$, and hence $U+m^{-1}\log|f|\le O(1)$ for any local equation $f$ of $mD$, by, for instance, the Siu decomposition of $dd^cU$.
\end{proof}

\begin{defi} A psh ray $U\colon\R_{>0}\to\PSH$ has \emph{algebraic singularities} if it induces a locally bounded psh metric on some normal, semiample test configuration $(\cX,\cL)$. 
\end{defi}
By Lemma~\ref{lem:pshtest}, such a ray $U$ has linear growth, and $U_\NA=\f_{(\cX,\cL)}\in\cHNA$. 

\begin{lem}\label{lem:algsing} For each $\f\in\cHNA$, there exists a smooth psh ray $U\colon\R_{\ge 0}\to\PSH$ with algebraic singularities such that $U_\NA=\f$. Further, every psh ray $V\colon\R_{>0}\to\PSH$ with $V_\NA\le\f$ satisfies $V\le U+O(1)$. 
\end{lem}
\begin{proof} By definition of $\cHNA$, we can pick a normal, semiample test configuration $(\cX,\cL)$ with $\f=\f_{(\cX,\cL)}$. Since $\cL$ is semiample, it admits a smooth $S^1$-invariant psh metric, which induces the desired psh ray $U$. If a psh ray $V$ satisfies $V_\NA\le\f$, then $V$ induces a psh metric on $(\cX,\cL)$ by Lemma~\ref{lem:pshtest}, and it follows that $V-U$ is bounded above. 
\end{proof}

%
%
%
%
\section{Ding-stability and twisted K\"ahler--Einstein currents}\label{sec:Dingcoerstab}
This section contains proofs of Theorems~A and~B in the introduction. In what follows, $(X,\om_0)$ is a compact K\"ahler manifold, $L$ an ample $\Q$-line bundle such that $\om_0\in c_1(L)$, and $\theta$ is a (quasi-positive) klt current with $c_1(X,\theta)=c_1(L)$. 
%
%
%
%
\subsection{Main results} 
The rest of this section will be devoted to the proof of the following result. 

\begin{thm}\label{thm:coerstab} If the Ding functional $\opD_\theta\colon\cE^1\to\R$ is coercive, then $(X,L)$ is uniformly Ding-stable with respect to $\theta$. If $\theta$ is further semipositive, the converse holds. 
\end{thm}

Combining this with Theorem~\ref{thm:Tboundcoer}, we obtain the following result, which is a more precise version of Theorem A in the introduction. 
\begin{cor}\label{cor:coerstab} If $\theta$ is semipositive, then: 
  \begin{itemize}
  \item[(i)] if $c_1(L)$ contains a $\theta$-twisted K\"ahler--Einstein current (resp.\ unique $\theta$-twisted K\"ahler--Einstein current), then $(X,L)$ is Ding-semistable (resp.\ uniformly Ding stable) with respect to $\theta$;  
  \item[(ii)] if $(X,L)$ is uniformly Ding-stable, then $c_1(L)$ contains a $\theta$-twisted K\"ahler--Einstein current.   
 \end{itemize}
 If we further assume that $\theta$ either has small unbounded locus, or is strictly positive, then the twisted K\"ahler--Einstein current in (ii) is further unique.
\end{cor}
\begin{proof}
  The statements in~(i) and~(ii) follow by combining Theorem~\ref{thm:Tboundcoer} and~\ref{thm:coerstab}. The uniqueness statement follows by also invoking Lemma~\ref{lem:uniqueKE}. 
\end{proof}
Theorem~B will be proved at the end of this section.
%
%
%
%
\subsection{Slopes of functionals}\label{sec:asympL}
Recall that each psh ray $U\colon\R_{>0}\to\PSH$ of linear growth induces a function $U_\NA\colon\Xdiv\to\R$, defined in terms of Lelong numbers. When $U$ has algebraic singularities, $U_\NA$ belongs to $\cHNA$, and we then have the following result, which is a reformulation of~\cite[Theorem 3.6]{BHJ2} (see also~\cite[Theorem 4.9]{Zak}, and~\cite{PRS} for a previous result in the same direction).

\begin{lem}\label{lem:slopeEJ} If $U\colon\R_{>0}\to\PSH$ is a psh ray with algebraic singularities, then 
\begin{itemize}
\item[(i)] $\opE(U_t)=t\ENA(U_\NA)+O(1)$; 
\item[(ii)]  $\opJ(U_t)=t\JNA(U_\NA)+O(1)$.
\end{itemize}
\end{lem}
Recall that we denote by the same letter a functional on $\cH$ and the induced functional on $\cHNA$
Coming back to the case of a general psh ray $U$, we set, as in~\eqref{eqn:DingNA},
$$
\LNA_\theta(U_\NA):=\inf_{\Xdiv}\left\{A_\theta+U_\NA\right\}\in\R\cup\{-\infty\}. 
$$
The following result is a generalization of~\cite[Proposition 3.8]{Berm16}, which basically corresponds to the case of algebraic singularities. The proof relies on a valuative criterion of integrability as discussed in Appendix~\ref{sec:valcrit}.
\begin{thm}\label{thm:Lasymp} For any psh ray $U\colon\R_{>0}\to\cE^1$ of linear growth, $\LNA_\theta(U_\NA)$ is finite, and coincides with the integrability threshold 
  $$
  \sup\left\{c\in\R\mid\int_1^{\infty}e^{2\left(c\,t-\opL_\theta(U_t)\right)}dt<\infty\right\}. 
  $$
\end{thm}
Further, when $\theta\ge 0$, the function $t\mapsto \opL_\theta(U_t)$ is convex (Lemma~\ref{lem:postwist}), and the integrability threshold is equal to its slope at infinity $\lim_{t\to\infty} t^{-1} \opL_\theta(U_t)$.

\smallskip
To prove Theorem~\ref{thm:Lasymp}, we may and do assume that $U$ extends to a quasi-psh function on $X\times\DD$, after adding, as before, a linear function of $t$. Then $U_\NA\le 0$, and hence $\LNA(U_\NA)\le\LNA(0)=0$, while the above integrability threshold is similarly nonpositive, since $L(U_t)\le O(1)$ for $t\gg 1$. 

In what follows, we denote for simplicity the log discrepancy function of $X\times\C$ by $A:=A_{X\times\C}$. Write $\theta=\theta_0+dd^c\p$ with $\theta_0$ smooth and $\p$ quasi-psh, and introduce the quasi-psh function
$$
V:=U+p_1^*\p
$$
on $X\times\DD$. Using $A_X(v)=A(\sigma(v))-1$ for $v\in\Xdiv$ and $\sigma$ the Gauss extension, we have
\begin{equation}\label{equ:DingW}
\LNA_\theta(U_\NA)=\inf_{w\in W}\left\{A(w)-w(V)\right\}-1
\end{equation}
with $W$ the set of all $\C^*$-invariant divisorial valuations $w$ on $X\times\C$ such that $w(\tau)=1$. 

\begin{lem}\label{lem:Lasymp} There exist $\e\in(0,1)$ and $C>0$ such that $w(V)\le(1-\e) A(w)+C$ for all $w\in W$. 
\end{lem}
\begin{proof} The restriction of the quasi-psh function $U$ on $X\times\DD$ to each submanifold $X\times\{\tau\}$ with $\tau\in\DD^*$ is in $\cE^1$; hence has zero Lelong numbers. Since Lelong numbers can only increase upon restriction, it follows that $U$ has zero Lelong number at each point of $X\times\DD^*$, and hence $e^{-U}\in L^q_{\loc}$ on $X\times\DD^*$ for every finite $q$, by Skoda's theorem. On the other hand, the assumption that $\theta$ is klt implies that $e^{-2\p}$ locally in $L^p_{\loc}$ for some $p>1$~\cite{Bernop,GuZh}. By H\"older's inequality, it follows that $e^{-2(1+\e)V}\in L^1_{\loc}$ on $X\times\DD^*$ for some $\e>0$. In other words, the multiplier ideal sheaf $\cJ((1+\e)V)$ is cosupported on $X\times\{0\}$, and hence contains some power of $\tau$, which yields $\sup_{w\in W} w(\cJ((1+\e)V))<\infty$. On the other hand, Lemma~\ref{lem:lesssing} shows that $w\left(\cJ((1+\e)V)\right)\ge (1+\e) w(V)-A(w)$ for all divisorial valuations $w$ on $X\times\C$ with $w(\tau)>0$, and the result follows. 
\end{proof}

\begin{proof}[Proof of Theorem~\ref{thm:Lasymp}] By definition of $\opL_\theta$, we have
$$
\opL_\theta(U_t)=-\tfrac 12\log\int_X e^{-2(V_t+\rho)}\om_0^n
$$
for some function $\rho\in C^\infty(X)$. Given $c\in\R$,  using $t=-\log|\tau|$, polar coordinates, and Fubini's theorem,  it is straightforward to see that the function $t\mapsto\exp\left(ct-\opL_\theta(U_t)\right)$ is $L^2$ in a neighborhood of $t=+\infty$ iff $|\tau|^{-c-1} e^{-V}$ is $L^2$ in a neighborhood of the central fiber in $X\times\DD$, or, equivalently, $L^2_{\mathrm{loc}}$ on $X\times\DD$. In view of~\eqref{equ:DingW}, we thus need to show
\begin{equation}\label{equ:Lasymp}
\inf\left\{s\in\R\mid |\tau|^s e^{-V}\in L^2_{\mathrm{loc}}\right\}=\sup_{w\in W}\left\{w(V)-A(w)\right\}<+\infty. 
\end{equation}

First suppose $|\tau|^se^{-V}\in L^2_{\mathrm{loc}}$. Applying Theorem~\ref{thm:valcritint1} (or Theorem~\ref{thm:valcritint2}) to $U=\log|\tau|$ shows that there exists $\e>0$ such that $s=w(U)\ge(1+\e)w(V)-A(w)\ge w(V)-A(w)$ for all $w\in W$.

For the reverse inequality  we use Theorem~\ref{thm:valcritint2} and Lemma~\ref{lem:Lasymp}.
Suppose $s=\sup_{w\in W}\{w(V)-A(w)\}+\d$, where $\d>0$. In particular, $s>-1$ as follows by taking $w$ as the order of vanishing along the central fiber.
We claim that $|\tau|^se^{-U}\in L^2_{\mathrm{loc}}$. By Theorem~\ref{thm:valcritint2}, it suffices to prove that there exists $\e'>0$ such that $s-(1+\e')w(V)+A(w)\ge 0$ for all $w\in W$.
Pick $\e\in(0,1)$ and $C\ge1$ as in Lemma~\ref{lem:Lasymp}.

If $A(w)\le 4C/\e$, then $w(V)\le(4/\e+1)C$ and hence 
\begin{align*}
  s-(1+\e')w(V)+A(w)
  &\ge\d+w(V)-A(w)-(1+\e')w(V)+A(w)\\
  &=\d-\e'w(V)\\
  &\ge\d-\e'(4/\e+1)C,
\end{align*}
which is nonnegative if $0<\e'\ll1$.
If instead $A(w)\ge 4C/\e$, then 
\begin{align*}
  s-(1+\e')w(V)+A(w)
  &>-1-(1+\e')(1-\e)(A(w)+C)+A(w)\\
  &\ge -1+\frac{\e}2A(w)-C,
\end{align*}
which again is nonnegative. This completes the proof.
\end{proof}
%
%
%
%
\subsection{Proof of Theorem~\ref{thm:coerstab}}\label{sec:proofDingRay}
Assume first that the Ding functional is coercive, \ie $\opD_\theta\ge\e J-C$ on $\cE^1$ for some $\e,C>0$. We then claim that $\DNA_\theta\ge\e\JNA$ on $\cHNA$, which will prove that $(X,L)$ is uniformly Ding-stable with respect to $\theta$. By Lemma~\ref{lem:algsing}, every $\f\in\cHNA$ is of the form $\f=U_\NA$ for some psh ray $U\colon\R_{>0}\to\cE^1$ with algebraic singularities. By Lemma~\ref{lem:slopeEJ}, we have $\opE(U_t)=t\ENA(U_\NA)+O(1)$ and $\opJ(U_t)=t\JNA(U_\NA)+O(1)$, while Theorem~\ref{thm:Lasymp} shows that $\LNA_\theta(\f)$ is the supremum of all $c\in\R$ such that $\int_1^\infty e^{2\left(c\,t-\opL_\theta(U_t)\right)}dt<\infty$. Now the coercivity assumption yields 
$$
\opL_\theta(U_t)\ge \opE(U_t)+\e \opJ(U_t)-C=t(\ENA(\f)+\e\JNA(\f))+O(1),
$$
and we infer $\LNA_\theta(\f)\ge\ENA(\f)+\e\JNA(\f)$; hence $\DNA_\theta\ge\e\JNA$ on $\cHNA$. 

\smallskip
Before proving the converse direction, let $U\colon\R_{>0}\to\cE^1$ be a psh ray with $U_t\le O(1)$ as $t\to\infty$, so that $U$ defines a quasi-psh function on $X\times\DD$ with multiplier ideals $\fa_m:=\cJ(m U)$ cosupported on the central fiber $X\times\{0\}$ (cf.~the proof of Lemma~\ref{lem:Lasymp}). By $S^1$-invariance of $U$, $\fa_m$ is $S^1$-invariant, and hence uniquely extends to a $\C^*$-invariant coherent ideal sheaf on $X\times\C$. 

\begin{lem}\label{lem:globgen} There exists $m_0\gg 1$ such that the sheaf
$\cO((m+m_0) p_1^*L)\otimes\fa_m$ is generated by its global
sections on $X\times\C$ for each $m\ge 1$. 
\end{lem}
\begin{proof} It is enough to show that $\cO((m+m_0) p_1^*L)\otimes\fa_m$ is $p_2$-globally generated, with $p_2\colon X\times\C\to\C$ denoting the second projection. We argue as in~\cite[Corollary 1.5]{DEL}. Pick a very ample line bundle $H$ on $X$, and choose $m_0$ such
that $A:=m_0L-K_X-(n+1)H$ is ample on $X$. By the relative version of the Castelnuovo-Mumford criterion, $\cO((m+m_0) p_1^*L)\otimes\fa_m$ is $p_2$-globally generated as soon as 
\begin{equation*}
  R^j (p_2)_*\left(\cO((m+m_0)p_1^*L-jp_1^*H)\otimes\fa_m\right)=0
\end{equation*}
for $1\le j\le n$, which holds away from $0\in\C$ by Kodaira vanishing, and near $0\in\C$ as a consequence of Nadel vanishing (compare~\cite[Theorem B.8]{siminag}). 
\end{proof}

\begin{lem}\label{lem:e201}
  Set $\f_m:=(m+m_0)^{-1}\f_{\fa_m}$. Then:  
\begin{itemize}
\item[(i)] $\f_m\in\cHNA$; 
\item[(ii)] $U_\NA\le\frac{m}{m+m_0}U_\NA\le\f_m\le\frac{m}{m+m_0}U_\NA+\frac{1}{m+m_0}(A_X+1)$ on $\Xdiv$; 
\item[(iii)] $\LNA_\theta(U_\NA)=\lim_{m\to+\infty}\LNA_\theta(\f_m)$. 
\end{itemize}
\end{lem}
\begin{proof} That $\f_m\in\cHNA$ is a direct consequence of Lemma~\ref{lem:globgen}. As for~(ii), the first inequality holds since $U_\NA\le 0$.  The remaining inequalities follow by applying the valuative analysis of multiplier ideals, see  Appendix~\ref{sec:valcrit}. Specifically, Lemma~\ref{lem:lesssing} gives 
  \begin{equation*}
  w(\cJ(mU))
  \le m\,w(U)
  \le w(\cJ(mU))+A(w),
\end{equation*}
for each divisorial valuation $w$ on $X\times\C$, and this implies the last two inequalities in~(ii) by setting $w=\sigma(v)$, where $v\in\Xdiv$ and $\sigma$ is the Gauss extension; indeed $A(w)=A_X(v)+1$.

To prove~(iii), first note that $\f_m\ge U_\NA$ implies
$$
\LNA_\theta(\f_m)
=\inf_{\Xdiv}\left\{A_\theta+\f_m\right\}
\ge\inf_{\Xdiv}\left\{A_\theta+U_\NA\right\}
=\LNA_\theta(U_\NA).
$$
To obtain an estimate in the opposite direction, let $\e>0$ and pick $v\in\Xdiv$ such that $A_\theta(v)+U_\NA(v)\le\LNA_\theta(U_\NA)+\e$. Then
$$
\LNA_\theta(U_\NA)\ge A_\theta(v)+U_\NA(v)-\e\ge\LNA(\f_m)+U_\NA(v)-\f_m(v)-\e, 
$$
which proves (iii) since $\f_m(v)\to U_\NA(v)$ by (ii). 
\end{proof}

\begin{lem}\label{lem:e202} For each $m$ we have $\ENA(\f_m)\ge\lim_{t\to+\infty}t^{-1} \opE(U_t)$. 
\end{lem}

\begin{rmk} Thus $\liminf\ENA(\f_m)\ge\lim_{t\to\infty}t^{-1}\opE(U_t)$. Here the inequality may be strict, see Example~\ref{exam:Eslope} below, based on~\cite{Dar17a}.
\end{rmk}

\begin{proof} By Lemma~\ref{lem:algsing}, we can choose a psh ray $U^m\colon\R_{>0}\to\cE^1$ with algebraic singularities such that $U^m_\NA=\f_m$, and hence $\opE(U^m_t)=t\ENA(\f_m)+O(1)$, in view of Lemma~\ref{lem:slopeEJ}. Since $U_\NA\le\f_m$, Lemma~\ref{lem:algsing} yields a constant $C>0$ such that 
$U_t\le U^m_t+C$ for $t\ge 1$. By monotonicity of $E$, we infer
$$
\opE(U_t)\le \opE(U^m_t)+O(1)=t\ENA(\f_m)+O(1),
$$ 
which concludes the proof
\end{proof}

We are now in a position to prove the reverse direction of Theorem~\ref{thm:coerstab}. Arguing by contradiction, assume that $\theta\ge 0$, $\DNA_\theta\ge\e\JNA$ on $\cHNA$ for some $\e\in(0,1)$, and that $\opD_\theta$ is not coercive. By Corollary~\ref{cor:coer}, we can then find a non-constant psh geodesic ray $U\colon\R_{\ge 0}\to\cE^1_{\sup}$ emanating from $0$ along which $\opM_\theta (U_t)\le 0$, and hence also $\opD_\theta(U_t)\le 0$, since $\opD_\theta\le \opM_\theta$. The assumptions on $U$ guarantee that $\opE(U_t)=ct$ for some $c<0$. As $D(U_t)=L(U_t)-\opE(U_t)\le 0$, we infer $L(U_t)\le ct$, and hence $\LNA(U_\NA)\le c$, by Theorem~\ref{thm:Lasymp}. Now consider the sequence $\f_m\in\cHNA$ constructed above. By Lemma~\ref{lem:slopesup}, we have $U_\NA(v_\triv)=0$; hence also $\sup\f_m=\f_m(v_\triv)=0$ by Lemma~\ref{lem:e201}~(ii). The assumption $\DNA_\theta\ge\e \JNA$ on $\cHNA$ thus yields 
$$
\LNA(\f_m)\ge(1-\e)\ENA(\f_m)
$$ 
for all $m$, and hence $\LNA(U_\NA)\ge(1-\e)c$, by Lemma~\ref{lem:e201} and Lemma~\ref{lem:e202}. We end up with $c\ge(1-\e)c$, a contradiction. 
%
%
%
%
\subsection{Proof of Theorem B}
Let $X$ be a projective manifold, $L$ an ample $\Q$-line bundle, and $\D$ an effective $\Q$-divisor with $(X,\D)$ klt and $c_1(L)=c_1(X,\D)$. By~\cite[Theorem 5.1]{BBEGZ}, the identity component of the algebraic group $\Aut^0(X,\D)$ acts transitively on $\D$-twisted K\"ahler--Einstein currents in $c_1(X,\D)$, and the stabilizer is compact by~\cite[Theorem 5.2]{BBEGZ}. If $c_1(X,\D)$ contains a unique such current $\om$, then $\Aut^0(X,\D)$ is contained in the compact group of isometries of $\om$, and is thus trivial, being an affine algebraic group. This proves (i)$\Longleftrightarrow$(ii). By Theorem~\ref{thm:DvsK}, $(X,\D)$ is uniformly K-stable iff it is uniformly Ding-stable with respect to $\D$. Thus (i)$\Longrightarrow$(iii) follows from (ii) in Corollary~\ref{cor:coerstab}. By Theorem~\ref{thm:coerstab}, (iii) conversely implies that the Ding functional $D_\D$ is coercive. Since the twisted Mabuchi functional $M_\D$ satisfies $M_\D\ge D_\D$, it is also coercive, and~\cite[Theorem 5.4]{BBEGZ} shows that $c_1(X,\D)$ contains a unique twisted K\"ahler--Einstein current;  hence (iii)$\Longrightarrow$(i). 

%
%
%
%
\section{Non-Archimedean potentials of finite energy and geodesic rays}\label{S202}
As above, $(X,\om_0)$ is a compact K\"ahler manifold, and $L$ is an ample $\Q$-line bundle with $\om_0\in c_1(L)$. Using part of the proof of Theorem~\ref{thm:coerstab}, we now undertake a deeper study the relationship between psh rays and non-Archimedean $L$-psh functions, and prove Theorem D in the introduction. 
%
%
\subsection{The Berkovich analytification}
Denote by $\XNA$ the \emph{Berkovich analytification}\footnote{This is usually denoted $\Xan$ in the literature~\cite{BerkBook}.} of $X$ with respect to the trivial absolute value on the ground field $\C$. We view $\XNA$ as a topological space, whose points can be understood as semivaluations on $X$, \ie valuations $v\colon\C(Y)^*\to\R$ on the function field of subvarieties $Y$ of $X$, trivial on $\C$. In particular, $\XNA$ contains the set $\Xdiv$ of divisorial valuations on $\C(X)$. Recall that, by convention, $\Xdiv$ contains the trivial valuation of $\C(X)$, denoted by $v_\triv$. The topology of $\XNA$ is generated by functions of the form $v\mapsto v(f)$ with $f$ a regular function on some Zariski open set $U\subset X$, and one shows that $\XNA$ is compact (Hausdorff), and that $\Xdiv\subset \XNA$ is dense. The projection $p_1\colon X\times\C\to X$ induces a map
$(X\times\C)^\NA\to\XNA$ that has a canonical continuous section, the Gauss extension
\begin{equation*}
  \sigma\colon\XNA\to(X\times\C)^\NA,
\end{equation*}
extending the map in~\S\ref{S302}. Its image
consists of all $\C^*$-invariant semivaluations $w$ satisfying $w(\tau)=1$.
%
%
\subsection{$L$-psh functions and psh rays}
As explained in~\cite{trivval}, any test configuration $(\cX,\cL)$ for $(X,L)$ defines a continuous metric on the Berkovich analytification of $L$. By subtracting the \emph{trivial metric}, defined by the trivial test configuration, we obtain a continuous function $\f_{(\cX,\cL)}\colon\XNA\to\R$ whose restriction to the dense subset $\Xdiv$ is the function defined in~\S\ref{S201}.

This allows us to view the elements of $\cHNA$ as continuous functions on all of $\XNA$. Concretely,
this can be explained as follows. Let $\fa$ be a $\C^*$-invariant ideal on $X\times\C$,
and write $\fa=\sum_{i\in\N}\tau^i\fa_i$ with $\fa_i$ ideals on $X$.
The function $\f_\fa\colon\XNA\to[-\infty,+\infty)$ given by 
$$
\f_\fa(v):=-\sigma(v)(\fa)=\max_i\{-v(\fa_i)-i\}. 
$$
is continuous, and finite-valued iff $\fa$ is cosupported on $X\times\{0\}$. 
This applies in particular to functions in $\cHNA$, which are of the form $\f=m^{-1}\f_\fa+c$ with $\fa$ cosupported on $X\times\{0\}$,  $p_1^*(mL)\otimes\fa$ globally generated, and $c\in\Q$.

An \emph{$L$-psh function} is a function $\f\colon\XNA\to[-\infty,+\infty)$, not identically $-\infty$, that can be written as the limit of a decreasing sequence in $\cHNA$. These functions are usc, satisfy the `maximum principle' 
\begin{equation}\label{equ:supLpsh}
\sup_{\XNA}\f=\f(v_\triv), 
\end{equation}
and are uniquely determined by their (finite) values on $\Xdiv$. The space $\PSHNA=\PSHNA(X,L)$ of $L$-psh functions is closed under decreasing limits. It is endowed with the weak topology of pointwise convergence on $\Xdiv$, and it is proved in~\cite{trivval}, as a consequence of~\cite{siminag}, that the space of sup-normalized functions
$$
\PSHNA_{\sup}:=\left\{\f\in\PSHNA\mid\sup\f=\f(v_\triv)=0\right\}
$$
is compact. 

\begin{lem}\label{lem:Lpsh} Let $m\ge 1$, and let $\fa$ be a $\C^*$-invariant coherent ideal sheaf on $X\times\C$ such that $mL$ is a line bundle and $p_1^*(mL)\otimes\fa$ is globally generated. Then $m^{-1}\f_\fa$ is $L$-psh. 
  \end{lem}
\begin{proof} For each $r\in\N$ we have $\f_{\fa^r}=r\f_{\fa}$. After replacing $m$ with a large enough multiple $rm$, we may thus assume that $mL$ is globally generated as well. For each integer $k\ge 1$, the $\C^*$-invariant ideal $\fa_k:=\fa+(\tau^k)$ is cosupported on $X\times\{0\}$, and $p_1^*(mL)\otimes\fa_k$ is globally generated since $p_1^*(mL)\otimes(\tau^k)$ and $p_1^*(mL)\otimes\fa$ are both globally generated. As a result, $m^{-1}\f_{\fa_k}\in\cHNA$, and we get the desired result since $\f_{\fa_k}=\max\{\f_{\fa},-k\}$ decreases pointwise to $\f_{\fa}$. 
\end{proof}

\begin{thm}\label{thm:rayLpsh}
  For each psh ray $U\colon\R_{>0}\to\PSH$ of linear growth, the function $U_\NA\colon\Xdiv\to\R$ 
  admits a unique extension to a function in $\PSHNA$.
\end{thm}
\begin{proof} Uniqueness follows from the fact that $L$-psh functions are determined by their restriction to $\Xdiv$. After adding to $U$ a linear function of $t$, we may as usual assume that it extends to a quasi-psh function on $X\times\DD$. By homogeneity, we may also assume that $L$ is an actual line bundle. 

For each $m\in\N$, the multiplier ideal sheaf $\fa_m:=\cJ(mU)$ can be viewed as a $\C^*$-invariant ideal sheaf on $X\times\C$, by $S^1$-invariance of $U$, and the proof of Lemma~\ref{lem:globgen} applies without change to yield $m_0\in\N$ such that $\cO((m+m_0) p_1^*L)\otimes\fa_m$ is globally generated for all $m$. As a result, 
$$
\f_m:=(m+m_0)^{-1}\f_{\fa_m}
$$ 
is $L$-psh, by Lemma~\ref{lem:Lpsh}. As in Lemma~\ref{lem:e201},  we further have 
$$
m U_\NA\le (m+m_0)\f_m\le m U_\NA+A_X+1
$$
on $\Xdiv$, which proves that $\f_m$ converges pointwise to $U_\NA$ on $\Xdiv$. Finally, the subadditivity property of multiplier ideals yields $\fa_{2m}\subset\fa_m^2$; hence 
$$
\f_{2m}\le\frac{2m+2m_0}{2m+m_0}\f_m\le\f_m,
$$
since $\f_m\le 0$. All in all, $\p_j:=\f_{2^j}$ is a decreasing sequence of $L$-psh functions, converging pointwise to $U_\NA$ on $\Xdiv$, and we conclude as desired that $U_\NA\in\PSHNA$. 
\end{proof}
%
%
\subsection{$L$-psh functions of finite energy}\label{sec:NAfin}
As in the complex case, the non-Archimedean Monge--Amp\`ere energy $\ENA\colon\cHNA\to\R$ defined in \S\ref{S204} admits a unique extension to a monotone, usc functional 
$$
\ENA\colon\PSHNA\to[-\infty,+\infty),
$$ 
obtained by setting for each $L$-psh function $\f$
$$
\ENA(\f)=\inf\left\{\ENA(\p)\mid\p\in\cHNA,\,\p\ge\f\right\}.
$$
We say that $\f$ has \emph{finite energy} if $\ENA(\f)>-\infty$ and write $\cENA$ for the space of such functions. To any $\f\in\cENA$ is attached a non-Archimedean Monge--Amp\`ere measure $\MA(\f)$, a Radon probability measure on $\XNA$. 

By the non-Archimedean Calabi--Yau theorem proved in~\cite{trivval} (building on~\cite{nama}), the non-Archimedean Monge--Amp\`ere operator sets up a one-to-one correspondence between $\cENA/\R$ and the set $\cMoneNA$ of Radon probability measures $\mu$ of finite energy, \ie such that 
\begin{equation}\label{eqn:ENAstar}
  \ENAstar(\mu):=\sup_{\f\in\cENA}(\ENA(\f)-\int\f\,d\mu)<\infty;
\end{equation}
given $\mu\in\cMoneNA$, the supremum in~\eqref{eqn:ENAstar} is attained for a unique $\f\in\cENA$, which then satisfies $\MA(\f)=\mu$. Conversely, we have 
\begin{equation}\label{eqn:ENAstardual}
  \ENA(\f)=\inf\{\ENAstar(\mu)+\int\f\,d\mu)\mid \mu\in\cMoneNA\},
\end{equation}
and the infimum is attained uniquely for $\mu=\MA(\f)$. As a consequence, we have
\begin{lem}\label{lem:ENA}
  Any two $\f,\p\in\cENA$ with $\f\ge\p$ satisfy $\ENA(\f)\ge\ENA(\p)$, with equality iff $\f=\p$. 
\end{lem}
\begin{proof}
  The inequality $\ENA(\f)\ge\ENA(\p)$ is clear from the definition. Now suppose $\f\ge\p$ and $\opE(\f)=\opE(\p)$. Set $\mu=\MA(\f)$. Then $\opE(\p)-\int\p\,\mu\ge \opE(\f)-\int\f\,\mu$, so since the supremum in~\eqref{eqn:ENAstar} is obtained uniquely for $\f$, we must have $\p=\f$.
\end{proof}

%
%
 \subsection{Maximal geodesic rays}\label{sec:maxray}
By Proposition~\ref{prop:georaylin}, any psh geodesic ray $U\colon\R_{\ge 0}\to\cE^1$ has linear growth; by Theorem~\ref{thm:rayLpsh}, it thus gives rise to an $L$-psh function $U_\NA\in\PSHNA$, and the following result implies that $U_\NA$ has finite energy. 

\begin{thm}\label{thm:slopeE} For any psh ray $U\colon\R_{>0}\to\cE^1$ of linear growth, the associated $L$-psh function $U_{\NA}$ belongs to $\cENA$, and 
\begin{equation}\label{equ:ineqslope}
\ENA(U_{\NA})\ge\lim_{t\to+\infty} t^{-1} \opE(U_t)>-\infty. 
\end{equation}
\end{thm}
The inequality can be strict in general, even for geodesic rays---see Example~\ref{exam:Eslope} below. 
\begin{proof} Using the notation of the proof of Theorem~\ref{thm:rayLpsh}, $\p_j:=\f_{2^j}$ is a decreasing sequence of functions in $\cHNA$, converging pointwise to $U_\NA$. By Lemma~\ref{lem:e202}, we further have, for each $j$, $\ENA(\p_j)\ge\lim_{t\to+\infty} t^{-1} \opE(U_t)$, which yields the desired result by continuity of $\ENA$ along decreasing sequences. 
\end{proof}

We now conversely show how to attach to each $\f\in\cENA$ a geodesic ray in $\cE^1$. 

\begin{defi} We say that a psh geodesic ray $U\colon\R_{\ge 0}\to\cE^1$ is \emph{maximal} if any psh ray of linear growth $V\colon\R_{>0}\to\cE^1$ with $\lim_{t\to 0} V_t\le U_0$ and $V_\NA\le U_\NA$ satisfies $V\le U$.
\end{defi}
A maximal geodesic ray is thus uniquely determined by $U_0$ and $U_\NA$. Not every psh geodesic ray is maximal, see Example~\ref{exam:Eslope} below.

\begin{thm}\label{thm:maxray} 
  For any $u\in\cE^1$ and any $\f\in\cENA$, there exists a unique maximal geodesic ray 
  $U\colon\R_{\ge 0}\to\cE^1$ emanating from $u$ such that $U_{\NA}=\f$. 
\end{thm}

\begin{proof} As already noticed, uniqueness is clear, so we need only prove existence.
First assume $u\in\cH$ and $\f\in\cHNA$. By Lemma~\ref{lem:algsing}, the set of smooth psh rays $V\colon\R_{\ge 0}\to\PSH$ with $V_0=u$ and $V_\NA=\f$ is non-empty; its usc upper envelope defines a psh geodesic ray $U\colon\R_{\ge 0}\to\cE^1$ with algebraic singularities such that $U_0=u$ and $U_\NA=\f$, by~\cite[Proposition 2.7]{Berm16}, and Lemma~\ref{lem:pshtest} shows that $U$ is maximal. 

Now consider the general case. Write $u$ and $\f$ as the limits of decreasing 
  sequences $u^j\in\cH$ and $\f^j\in\cHNA$, respectively.
  For each $j$, we have a maximal geodesic ray $U^j$ with 
  $U^j_0=u_j$ and $U^j_\NA=\f_j$. By maximality, $U^{j+1}\le U^j$, so the limit $U:=\lim_jU^j$
  exists. By Corollary~\ref{cor:geod} and each Lemma~\ref{lem:slopeEJ}, we have, for each $j$ and $t$,
  $$
  \opE(U^j_t)=\opE(u^j)+t\ENA(\f^j)\ge \opE(u)+t\ENA(\f)>-\infty,
  $$ 
  so $U_t\in\cE^1$ and $\opE(U_t)=\opE(u)+t\ENA(\f)$.
  Thus $U\colon\R_{\ge 0}\to\cE^1$ is a psh geodesic ray, by Corollary~\ref{cor:geod}. 
  On the one hand, $U\le U^j$ implies $U_\NA\le U^j_\NA=\f_j$ for all $j$, and 
  hence $U_\NA\le\f$. 
  On the other hand, the formula $\opE(U_t)=\opE(u)+t\ENA(\f)$ yields 
  $\ENA(U_\NA)\ge\ENA(\f)$ by Theorem~\ref{thm:slopeE}, and hence $U_\NA=\f$, by Lemma~\ref{lem:ENA}. 

  Finally, suppose $V\colon\R_{> 0}\to\cE^1$ is a psh ray of linear growth with $\lim_{t\to 0}V_t\le u$
  and $V_\NA\le\f$. Since $u\le u_j$ and $U_\NA\le\f_j$, we have $V\le U^j$ by 
  maximality of $U^j$, and hence $V\le U$.
\end{proof}

\begin{cor}\label{cor:maxraychar}
  A psh geodesic ray $U\colon\R_{\ge 0}\to\cE^1$ is maximal iff equality holds in~\eqref{equ:ineqslope}, or, equivalently, $\opE(U_t)=\opE(U_0)+t\ENA(U_\NA)$ for all $t\ge0$.
\end{cor}
\begin{proof} Since $\opE(U_t)$ is an affine function of $t$, $\lim_{t\to\infty} t^{-1} \opE(U_t)=\ENA(U_\NA)$ is equivalent to $\opE(U_t)=\opE(U_0)+t\ENA(U_\NA)$, and the proof of Theorem~\ref{thm:maxray} shows that 
  the latter holds when $U$ is maximal. Assume, conversely, that $\opE(U_t)=\opE(U_0)+t\ENA(U_\NA)$ for all $t$, and let $U'$ be the maximal geodesic ray with $U'_0=U_0$ and $U'_\NA=U_\NA$. 
  Then $U\le U'$, and, as we have seen, $\opE(U'_t)=\opE(U_0)+t\ENA(U_\NA)$ for all $t$. 
  For each $t\ge 0$, we thus have $U_t\le U'_t$ and $\opE(U_t)=\opE(U'_t)$, which yields $U_t=U'_t$, proving that $U=U'$ is maximal.  
\end{proof}

\begin{exam}\label{exam:algmax} By Lemma~\ref{lem:slopeEJ}, every psh geodesic ray $U$ with algebraic singularities is maximal. Conversely, a maximal geodesic ray $U$ has algebraic singularities iff $U_\NA$ belongs to $\cHNA$. 
\end{exam}

\begin{exam} By~\cite{nakstab}, every linearly bounded filtration $\cF$ of the algebra of sections 
$$
R(X,L)=\bigoplus_{m\in\N} H^0(X,mL)
$$ 
gives rise to a bounded $L$-psh function $\f$ on $\XNA$. On the other hand, Ross and Witt Nystr\"om associate to $\cF$ a psh geodesic ray $U$ emanating from $0$~\cite[Corollary 7.12]{RWN}, and one can check that $U$ is indeed the maximal geodesic ray with $U_\NA=\f$. 
\end{exam}

\begin{exam}\label{exam:Eslope} Let $X=\P^1$, $L=\cO(1)$ and $\om\in c_1(L)$ the Fubini-Study metric. Following~\eg~\cite[Thm~3, p.31]{Car67}, we can construct a polar Cantor set $K\subset\P^1$. This carries an atom-free probability measure $\mu$, whose   potential $v\in\PSH(X,\om)$ has no Lelong numbers (because $\mu$ has no atoms), but does not belong to the class $\cE$ defined in~\cite{GZ} (since $\mu$ has positive mass on the polar set $K$). Now use $v$ to construct a psh geodesic ray $U$ emanating from $0$ as in~\cite[Theorem~2]{Dar17a}. Since $v$ has zero Lelong numbers, so does $U$, so $U_\NA=0$. However, $U$ is not constant by~\cite[Theorem~4.1]{Dar17a}, and hence not maximal by Corollary~\ref{cor:maxraychar}.
\end{exam}
%
%
%
\subsection{Uniform Ding-stability, reprise}
Ding-stability of $(X,L)$ with respect to $\theta$ was defined in \S\ref{S204} in terms of the non-Archimedean Ding functional $\DNA_\theta$ on $\cHNA$. As in~\cite[Lemma 2.9]{trivval}, we first show that it can equivalently be formulated as a condition on the whole space $\cENA$. 

\begin{lem}\label{lem:NADingE1} Given any klt current $\theta$, $(X,L)$ is Ding-semistable (resp.~uniformly Ding-stable) with respect to $\theta$ iff $\DNA_\theta\ge 0$ on $\cENA$ (resp.\ $\DNA_\theta\ge\e\JNA$ on $\cENA$ for some $\e>0$). 
\end{lem}
\begin{proof} Given $\e\ge 0$ and $\f\in\cENA$, $\DNA_\theta(\f)\ge\e\JNA(\f)$ is equivalent to 
$$
A_\theta(v)+\f(v)\ge(1-\e)\ENA(\f)+\e\f(v_\triv)
$$
for all $v\in\Xdiv$. If this holds for all $\f\in\cHNA$, then it also holds for $\f\in\cENA$, by continuity of $\ENA$ along decreasing sequences.  
\end{proof}

Using the results of \S\ref{sec:maxray}, we are now in a position to prove Theorem D, which we reformulate here for convenience.
\begin{thm}\label{thm:unifDing} Let $\theta$ be a semipositive klt current such that $c_1(X,\theta)=c_1(L)$. The following are equivalent: 
\begin{itemize}
\item[(i)] $\opD_\theta\colon\cE^1\to\R$ is coercive;
\item[(ii)] $\DNA_\theta(\f)>0$ for all non-constant $\f\in\cENA$;
\item[(iii)] $(X,L)$ is uniformly Ding-stable with respect to $\theta$.
\end{itemize}
\end{thm}
\begin{proof} (i)$\Longleftrightarrow$(iii) is the content of Theorem~\ref{thm:coerstab}, and Lemma~\ref{lem:NADingE1} shows that (ii)$\Longrightarrow$(iii). Now assume (ii), and suppose by contradiction that (i) fails. By Theorem~\ref{thm:coer}, 
$$
\opL_\theta(U_t)-\opE(U_t)=\opD_\theta(U_t)\le 0
$$ 
for some non-constant psh geodesic ray $U\colon\R_{\ge 0}\to\cE^1_{\sup}$, which thus satisfies $\opE(U_t)=c t$ for all $t$, where $c<0$. By Theorem~\ref{thm:Lasymp} and Theorem~\ref{thm:slopeE}, we infer
$$
\LNA_\theta(U_\NA)\le c\le\ENA(U_\NA).
$$
Thus $U_\NA\in\cENA$ satisfies $\DNA_\theta(U_\NA)=\LNA(U_\NA)-\ENA(U_\NA)\le 0$ and is sup-normalized, see Lemma~\ref{lem:slopesup}, so~(ii) yields $U_\NA=0$, which contradicts $\LNA_\theta(U_\NA)\le c<0$.
\end{proof}
%
%
%
%
\section{The stability threshold and the greatest Ricci lower bound}
As before, $X$ is a smooth projective variety with an ample $\Q$-line bundle $L$. Following~\cite{FO,BlJ,nakstab}, we characterize Ding-stability with respect to a klt current in terms of a stability threshold, and then prove Theorem C. 
%
%
%
\subsection{The expected vanishing order}
Assume first that $L$ is an actual line bundle (as opposed to a $\Q$-line bundle). Given a valuation $v\in\Xdiv$ and a nonzero section $s\in H^0(X,L)$, we can make sense of $v(s)\in\Q_{\ge 0}$, by evaluating $v$ on the local function corresponding to $s$ in a trivialization of $L$ at the center of $v$. This defines a filtration $F^\la:=\{s\mid v(s)\ge\la\}$ of $H^0(X,L)$, which shows that $v$ takes only finitely many values $\la\in\Q_{\ge 0}$ on $H^0(X,L)\setminus\{0\}$, and provides a way to count these with multiplicity $\dim\mathrm{Gr}^\la$, giving rise to the \emph{vanishing sequence} of $v$ on $L$~\cite{BKMS}. 

For each $m\in\N$, define $S_m(L)$ as the mean value of the vanishing sequence of $mL$, divided by $m$. By~\cite[Lemma 3.5]{BlJ}, we have
\begin{equation}\label{equ:Sm}
S_m(v)=\max\left\{v(D)\mid D\text{ of $m$-basis type}\right\},
\end{equation}
where a divisor of \emph{$m$-basis type for $L$} is a $\Q$-divisor of the form 
$$
D=\frac{1}{mN_m}\sum_{j=1}^{N_m}\div(s_j)
$$ 
for some basis $(s_1,\dots,s_{N_m})$ of $H^0(mL)$. By~\cite{BC,BKMS}, the vanishing sequence of $v$ on $mL$, scaled by $1/m$, equidistributes as $m\to\infty$. The sequence $S_m(v)$ thus admits a limit $S_L(v)\in\R_{>0}$, the \emph{expected vanishing order} of multisections of $L$ along $v$. 

By~\cite[\S 2.4]{BKMS} and~\cite[Lemma 5.13]{BHJ1}, this invariant can be expressed as 
\begin{equation}\label{equ:SL}
S_L(v)=V^{-1}\int_0^{+\infty}\vol\left(L,v\ge\la\right)\,d\la,
\end{equation}
where $\vol(L,v\ge\la)$ denotes the volume of the graded subalgebra of the section ring $R(X,L)=\bigoplus_{m=0}^\infty H^0(X,mL)$ consisting of sections $s\in H^0(X,mL)$ such that $v(s)\ge m\la$ (see~\textit{loc.\,cit.} for details). In particular, if $x\in X$, then $S_L(\ord_x)$ coincides with the invariant considered in~\cite[\S 4]{MR}. 

By construction, $S_L(v)$ is homogeneous of degree $1$ with respect to $L$, and it can thus be defined for $L$ a $\Q$-line bundle, by setting $S_L(v):=m^{-1} S_{mL}(v)$ for any $m\in\Z_{>0}$ such that $mL$ is a line bundle. 

A key point for what follows is that the convergence of $S_m(v)$ of $S_L(v)$ is actually semiuniform, in the following sense: 
\begin{lem}\label{lem:convS}\cite[Corollary 3.6]{BlJ}
For each $\e>0$, there exists $m_0$ such that $S_m(v)\le(1+\e)S_L(v)$ for all $m\ge m_0$ and all $v\in\Xdiv$. 
\end{lem} 
%
%
%
\subsection{The stability threshold}\label{sec:stabthre}

Following~\cite{FO,BlJ,nakstab}, we introduce: 

\begin{defi}\label{defi:stab} Given a klt current $\theta$, we define the \emph{stability threshold} of $(X,L)$ with respect to $\theta$ as
$$
\d_\theta(X,L):=\inf_{v\in\Xdiv}\frac{A_\theta(v)}{S_L(v)}
$$
\end{defi}
When $\theta=0$, we simply write $\d(X,L)$, and recover the invariant studied in~\cite{BlJ,nakstab}. Since the latter is positive, so is $\d_\theta(X,L)$, as follows from Corollary~\ref{cor:lctval1}. Note also that $\d_\theta(X,tL)=t^{-1}\d_\theta(X,L)$ for $t\in\Q_{>0}$. 

On the other hand, the \emph{log canonical threshold} of an effective $\Q$-divisor $D$ with respect to $\theta$ is defined as 
\begin{equation}\label{equ:lct}
\lct_\theta(D):=\sup\left\{c\ge 0\mid\cJ\left(\theta+c\,D\right)=\cO_X\right\}=\inf_{v_\triv\ne v\in\Xdiv}\frac{A_\theta(v)}{v(D)};
\end{equation}
see Corollary~\ref{cor:lctval2} for the second equality. Adapting, respectively, the arguments of~\cite[Theorem 4.4]{BlJ} and~\cite[Theorem 2.14]{nakstab}, we will prove: 

\begin{thm}\label{thm:stab} The twisted stability threshold satisfies the following properties: 
\begin{itemize}
\item[(i)] $\d_\theta(X,L)$ is the limit as $m\to\infty$ of 
$$
\d^{(m)}_\theta(X,L):=\inf\left\{\lct_\theta(D)\mid D\text{ of $m$-basis type }\right\};
$$
  \item[(ii)] $(X,L)$ is Ding-semistable (resp.\ uniformly Ding-stable) with respect to $\theta$ iff $\d_\theta(X,L)\ge 1$ (resp.\ $\d_\theta(X,L)>1$). 
\end{itemize}
\end{thm}
When $\theta=0$, (ii) follows from~\cite[Theorem 1.3]{Fuj19} (note that the invariant $\b(v)$ therein is $A_X(v)-S_L(v)$ multiplied by $V$, by~\eqref{equ:SL}). While Fujita's arguments rely on the Minimal Model Program, our proof of Theorem~\ref{thm:stab} builds on the non-Archimedean analogue of the thermodynamical formalism (compare Lemma~\ref{lem:legendre}), as in~\cite{nakstab}.
\begin{lem}\label{lem:Dirac}
  For each $v\in\Xdiv$, there exists a unique $\f_v\in\cENA$ such that $\MA(\f_v)=\d_v$ and $\f_v(v)=0$. Further, we have $\ENA(\f_v)=\ENAstar(\d_v)=S_L(v)$ and $\opJ(f_v)\ge n^{-1}S_L(v)$.
\end{lem}
\begin{proof}
  By~\cite[Proposition~5.6]{nakstab}, any valuation $v\in\Xdiv$ is nonpluripolar in the sense that $\f(v)>-\infty$ for all $\f\in\PSH$. The existence of $\f_v$ is therefore a special case of~\cite[Theorem~5.13]{nakstab}, which also gives the formulas $\ENA(\f_v)=\ENAstar(\d_v)=S_L(v)$ and $\JNA(\f_v)=T_L(v)-S_L(v)$, where $T_L(v)=-\inf\{\f(v)\mid \in\PSH_{\sup}\}<\infty$. The inequality $\JNA(\f_v)\ge n^{-1}S_L(v)$ now follows from equation~(5.3) in~\cite{nakstab}.
\end{proof}

\begin{proof}[Proof of Theorem~\ref{thm:stab}] By~\eqref{equ:Sm} and~\eqref{equ:lct}, we have 
$\d^{(m)}_\theta(X,L)=\inf_{v\in\Xdiv}\frac{A_\theta(v)}{S_m(v)}$, and hence 
$$
\limsup\d^{(m)}_\theta(X,L)\le\d_\theta(X,L).
$$
On the other hand, for each $\e>0$ we have $S_m\le(1+\e) S_L$ on $\Xdiv$ for all $m\gg 1$, thanks to Lemma~\ref{lem:convS}. This implies $\d^{(m)}_\theta(X,L)\ge(1+\e)^{-1}\d_\theta(X,L)$, and proves (i). 

To prove (ii), assume first $(X,L)$ Ding-semistable (resp.\ uniformly Ding-stable) with respect to $\theta$. By Lemma~\ref{lem:NADingE1}, $\DNA_\theta\ge\e\JNA$ on $\cENA$ with $\e\ge 0$ (resp.\ $\e>0$). Pick any $v\in\Xdiv$, and consider $\f_v\in\cENA$ as in Lemma~\ref{lem:Dirac}. The inequality $\DNA_\theta(\f_v)\ge\e\JNA(\f_v)$ gives
\begin{equation*}
  A_\theta(v)-S_L(v)
  \ge\inf_{\Xdiv}(A_\theta+\f_v)-S_L(v)
  =\opD_\theta(\f_v)
  \ge\e\opJ(f_v)
  \ge\e n^{-1}S_L(v),
\end{equation*}
so $A_\theta(v)\ge(1+\e n^{-1})S_L(v)$ for all $v\in\Xdiv$, and hence $\d_\theta(X,L)\ge 1+\e n^{-1}$. 

Conversely, assume $\d_\theta(X,L)\ge\d$ for some $\d\in\Q\cap[1,+\infty)$, \ie $A_\theta(v)\ge\d S(v)$ for $v\in\Xdiv$, and pick $\f\in\cHNA$ with $\sup\f=0$. Then $\d^{-1}\f\in\cHNA$, so~\eqref{eqn:ENAstardual} yields
\begin{multline*}
  \ENA(\d^{-1}\f)
  =\inf_{\mu\in\cMoneNA}(\opE^*(\mu)+\int\d^{-1}\f\,\mu)
  \le\inf_{v\in\Xdiv}(\opE^*(\d_v)+\d^{-1}\f(v))\\
  =\inf_{v\in\Xdiv}(S_L(v)+\d^{-1}\f(v))
  \le\d^{-1}\inf_{v\in\Xdiv}(A_\theta(v)+\f(v))
  =\d^{-1}\opL_\theta(\f)
\end{multline*}
Combining this with the inequality $\d\ENA(\d^{-1}\f)\ge\d^{-1/n}\ENA(\f)$ from~\cite[Lemma 6.17]{trivval}, we get $\LNA_\theta(\f)\ge\d^{-1/n}\ENA(\f)$. Since $\sup\f=0$, we have $\opJ(\f)=-\opE(\f)$, and so
\begin{equation*}
  \DNA_\theta(\f)
  =\LNA_\theta(\f)-\ENA(\f)
  \ge(\d^{-1/n}-1)\ENA(\f)
  =(1-\d^{-1/n})\JNA(\f),
\end{equation*}
which completes the proof.
\end{proof}
%
%
\subsection{The greatest twisted Ricci lower bound}
In order to state the next result, we introduce the following invariants: 
\begin{itemize}
\item[(a)] the \emph{greatest twisted Ricci lower bound}
$$
\b_\theta(X,L):=\sup\left\{\b\in\R\mid\exists\om\in c_1(L),\,\Ric_\theta(\om)\ge\b\om\right\}; 
$$
\item[(b)] the \emph{nef threshold}
$$
s_\theta(X,L):=\max\left\{s\in\R\mid c_1(X,\theta)\ge sc_1(L)\right\}. 
$$
\end{itemize}
In (a), $\om$ is a current of finite energy in $c_1(L)$, and $\Ric(\om)\ge\b\om+\theta$ means that the difference is a smooth semipositive $(1,1)$-form. In (b), $c_1(X,\theta)\ge s c_1(L)$ means that the difference is nef. 

The next result is Theorem C in the introduction. 

\begin{thm}\label{thm:Riccistab} For any semipositive klt current $\theta$, we have 
$$
\b_\theta(X,L)=\min\{\d_\theta(X,L),s_\theta(X,L)\}. 
$$
\end{thm}
Note that we do not require $c_1(X,\theta)=c_1(L)$. In the usual Fano case $\theta=0$, $L=-K_X$, the nef threshold is clearly equal to $1$, and hence: 

\begin{cor}\label{cor:Riccistab} If $X$ is a Fano manifold $X$, then $\b(X)=\min\{\d(X),1\}$. In particular, $X$ is K-semistable iff for each K\"ahler form  $\om\in c_1(X)$ and $t\in(0,1)$ there exists a K\"ahler form $\om_t\in c_1(X)$ such that 
$$
\Ric(\om_t)=t\om_t+(1-t)\om.
$$
\end{cor}
This corollary was independently established in the appendix of~\cite{CRZ}, as a consequence of~\cite{LS,SW} (see also~\cite{Li11} for the toric case and~\cite{Cab19} for the case of Fano $\theta$-manifolds of complexity one). The final statement was also previously obtained in~\cite{Li17a}, also building on~\cite{CDS}.

\begin{proof}[Proof of Theorem~\ref{thm:Riccistab}] We obviously have $\b_\theta(X,L)\le s_\theta(X,L)$. Consider first $s>0$ with $c_1(X,\theta)+sc_1(L)$ ample, and pick a K\"ahler form $\a$ in this class. The equation $\Ric(\om_u)=-s\om_u+\theta+\a$ with $u\in\cE^1$ corresponds to a Monge--Amp\`ere equation of the form $\MA(u)=e^{2(su-\p-\rho)}\om_0^n$ with $\theta-dd^c\p$ and $\rho$ smooth, and hence admits a solution~\cite{BBGZ}. It follows that $s_\theta(X,L)\le 0\Longrightarrow s_\theta(X,L)=\b_\theta(X,L)$, which proves the theorem in that case. 

Assume now $s_\theta(X,L)>0$, and pick $s\in\Q_{>0}$ with $c_1(X,\theta)-s c_1(L)$ ample. If $\Ric(\om)=s\om+\theta+\a$ for some $\om\in c_1(L)$ and $\a\ge 0$, Corollary~\ref{cor:coerstab} shows that $(X,sL)$ is Ding-semistable with respect to $\theta+\a$, and hence with respect to $\theta$ as well, which yields $s\le s_\theta(X,L)$, and hence $\b_\theta(X,L)\le\min\{s_\theta(X,L),\d_\theta(X,L)\}$. 

Conversely, pick $s\in\Q_{>0}$ with $c_1(X,\theta)-s c_1(L)$ ample and $s<\d_\theta(X,L)$, \ie $(X,sL)$ uniformly Ding-stable with respect to $\theta$. For any choice of K\"ahler form $\a\in c_1(X,\theta)-s c_1(L)$, we have $c_1(X,\theta+\a)=c_1(sL)$, and Corollary~\ref{cor:coerstab} thus yields $\om\in c_1(L)$ solving $\Ric(\om)=s\om+\theta+\a$, which proves $\b_\theta(X,L)\ge\min\{s_\theta(X,L),\d_\theta(X,L)\}$. 
\end{proof}
%
%
%
%
\appendix
\section{Estimates}
In what follows, $C_n$ denotes a constant that only depends on the dimension $n=\dim X$, but whose value may change from line to line.
\begin{lem}\label{L301}
  If $u_j,v_j\in\cE^1$, $0\le j\le n$, then 
  \begin{equation*}
    \big|
    \int(u_0-v_0)(\om_{u_1}\wedge\dots\wedge\om_{u_n}-\om_{v_1}\wedge\dots\wedge\om_{v_n})
    \big|
    \le C_n
    \opI(u_0,v_0)^{\frac1{2^n}}
    \max_{1\le p\le n}\opI(u_p,v_p)^{\frac1{2^n}}
    M^{1-\frac1{2^{n-1}}},
  \end{equation*}
  where $M=\max_{0\le j\le n}\max\{\opI(u_j),\opI(v_j)\}$. 
\end{lem}
Here $\opI(u,v)=V^{-1}\int_X(u-v)(\om_v^n-\om_u^n)$ for $u,v\in\cE^1$. We also write $\opI(u):=\opI(u,0)$. 
\begin{proof}
  After regularization~\cite{BK07}, we may assume that $u_j,v_j\in\cH$ for all $j$. 
  For $0\le p\le n$, set 
  $\eta_p:=\om_{u_1}\wedge\dots\wedge\om_{u_p}\wedge\om_{v_{p+1}}\wedge\dots\wedge\om_{v_n}$ and $A_p:=\int(u_0-v_0)\eta_p$. Then we want to estimate $|A_n-A_0|$.
  Now $A_p-A_{p-1}=\int(u_0-v_0)dd^c(u_p-v_p)\wedge\eta$,
  so by Stokes and Cauchy--Schwartz, we have $|A_p-A_{p-1}|^2\le b_pc_p$, where 
  $b_p=\int d(u_0-v_0)\wedge d^c(u_0-v_0)\wedge\eta_p$
  and $c_p=\int d(u_p-v_p)\wedge d^c(u_p-v_p)\wedge\eta_p$.
  Set $w_p:=\frac1{n-1}(u_1+\dots+u_{p-1}+v_{p+1}+\dots+v_n)$.
  Then 
  \begin{equation*}
    b_p
    \le C_n\int d(u_0-v_0)\wedge d^c(u_0-v)\om_{w_p}^n
    \le C_n \opI(u_0,v_0)^{1-\frac1{2^{n-1}}}\max\{\opI(u_0,w_p),\opI(v_0,w_p)\}^{1-\frac1{2^{n-1}}},
  \end{equation*}
  where the second equality follows from~\cite[Lemma~1.9]{BBEGZ}.
  Now $\opI(u_0,w_p)\le C_n\max\{\opI(u_0),\opI(w_p)\}$. Since the $I$ and $J$ functionals are 
  comparable, and $u\mapsto \opJ(u)$ is convex, it easily follows that $\opI(w_p)\le C_n M$.
  Applying the analogous estimate with $v_0$ instead of $u_0$, we get 
  $b_p\le C_n\opI(u_0,v_0)^{\frac{1}{2^{n-1}}}M^{1-\frac1{2^{n-1}}}$.
  Similarly, $c_p\le C_n\opI(u_p,v_p)^{\frac{1}{2^{n-1}}}M^{1-\frac1{2^{n-1}}}$,
  and the result follows.
\end{proof}

\begin{lem}\label{L302}
  If $u_0,v_0\in\cE^1$ and $f_j\in\cE^1$, $0\le j\le n$, then 
  \begin{equation*}
    \big|
    \int(u_0-v_0)\om_{f_1}\wedge\dots\wedge\om_{f_n}
    \big|
    \le C_n
    d_1(u_0,v_0)^{\frac1{2^n}}M^{1-\frac1{2^n}},
  \end{equation*}
  where $M=\max\{\opI(u_0),\opI(v_0),\max_{1\le j\le n}\opI(f_j)\}$.
\end{lem}

\begin{proof}
  By Lemma~\ref{L301} we have 
  \begin{equation}\label{e301}
    |\int(u_0-v_0)(\om_{f_1}\wedge\dots\wedge\om_{f_n}-\om_{u_0}^n)|
    \le C_n\opI(u_0,v_0)^{\frac1{2^n}}M^{1-\frac1{2^n}}.
  \end{equation}
  Now~\cite[Theorem~5.5]{Dar15} shows that 
  \begin{equation*}
    C_n^{-1}d_1(u_0,v_0)
    \le
    \int|u_0-v_0|(\om_{u_0}^n+\om_{v_0}^n)
    \le C_nd_1(u_0,v_0).
  \end{equation*}
  This first implies that 
  $\opI(u_0,v_0)=\int(u_0-v_0)(\om_{u_0}^n+\om_{v_0}^n)\le C_nd_1(u_0,v_0)$,
  and then that 
  \begin{equation*}
    |\int(u_0-v_0)\om_{u_0}^n|
    \le
    C_n \opI(u_0,v_0)
    \le C_n \opI(u_0,v_0)^{\frac1{2^n}}\max\{\opI(u_0),\opI(v_0)\}^{1-\frac1{2^n}}
    \le C_n d_1(u_0,v_0)^{\frac1{2^n}}M^{1-\frac1{2^n}}.
  \end{equation*}
  Combining this with~\eqref{e301} completes the proof.
\end{proof}

\begin{cor}\label{cor:intd} If $u\in\cE^1$ then $\left|\int_X u\,\om^n\right|\le C_n d_1(u,0)$. 
\end{cor}
\begin{proof} This is a consequence of Lemma~\ref{L302}, since $\opI(u)\le C_n d_1(u,0)$. 
\end{proof}
%
%
%
%
\section{A valuative criterion of integrability}\label{sec:valcrit}
The study of the asymptotics of the Ding functional along psh rays in~\S\ref{sec:Dingcoerstab} relies on the valuative analysis of singularites of psh functions, as developed in~\cite{hiro} following earlier work in~\cite{pshsing,valmul} in dimension two. Here we revisit some of this analysis, taking advantage of the solution of the openness conjecture that was unknown at the time of~\textit{loc.\,cit.}
%
%
\subsection{Preliminaries}
Throughout this section, $X$ is a connected complex manifold of dimension $n$.
Recall that a \emph{quasi-psh} function on $X$ is a function that is locally the sum of a smooth function and a psh function.

A holomorphic map between complex manifolds is called a \emph{modification} if it is proper and bimeromorphic. For example, the analytification of a proper birational morphism between smooth algebraic varieties is a modification.
If $\mu\colon X'\to X$ is a modification (with $X'$ a complex manifold), we denote by $K_{X'/X}$ the relative canonical divisor, defined locally by the Jacobian determinant of $\mu$.

A modification $\mu\colon X'\to X$ is called \emph{projective} if for any point $x\in X$, there exists $N\ge 1$ and an open neighborhood $U$ of $x$ in $X$ such that the restriction of $\mu$ to $U'=\mu^{-1}(U)$ is the composition of a closed embedding $U'\to U\times\P^N$ and the projection to the first factor.

If $Z\subset X$ is any subset,  $\Ltwoloc(X,Z)$ denotes the set of Lebesgue measurable functions $h$ defined in a neighborhood of $Z$ in $X$ such that $|h|^2$ is locally
integrable in a neighborhood of every point $x\in Z$.
%
%
\subsection{Singularity classes of quasi-psh functions}
If $D$ is an effective divisor on $X$, then we say that a quasi-psh function $U$ on $X$ has \emph{divisorial singularities of type $D$} if for every $x\in X$, we have that $U=\log|f|+O(1)$ near $x$, where $f$ is a local equation for $D$. A partition of unity argument shows that given any effective divisor $D$ there exists a quasi-psh function $U$ with divisorial singularities of type $D$, and any two such functions differ by a locally bounded function, see~\eg~\cite[Prop.\ 1.1.4]{SebThesis}.

An effective $\R$-divisor on $X$ is a formal sum $D=\sum_{i=1}^mc_iD_i$, where $D_i$ is an effective divisor and $c_i\in\R_{\ge0}$ for all $i$. We say that a quasi-psh function $U$ has divisorial singularities of type $D$ if, locally, $U=\sum c_i\log|f_i|+O(1)$, where $f_i$ is a local equation for $D_i$. If $D'$ is an effective divisor whose support contains the $D_i$, then we also say that $U$ has \emph{divisorial singularities along $D'$}.

More generally, consider coherent ideal sheaves $\fa_1,\dots,\fa_m$ and real numbers $c_1,\dots,c_m\ge0$. We say that a quasi-psh function $U$ has \emph{analytic singularities of type $\prod_i\fa_i^{c_i}$} if locally, $U=\sum_ic_i\log|\fa_i|+O(1)$. Here we write $\log|\fa|$ for the function defined locally near $x\in X$ by $\max_j\log|f_j|$, where $\{f_j\}_j$ is a set of generators of $\fa_x$; this is well defined up to $O(1)$. The class of quasi-psh functions with analytic singularities is stable under finite sums.

It follows from Hironaka's theorem that given any quasi-psh function $U$ on $X$ with analytic singularities, there exists a projective modification $\mu\colon X'\to X$ such that $U\circ\mu$ has divisorial singularities along an snc divisor whose support contains the support of $K_{X'/X}$. We call $\mu$ a \emph{log resolution} of the singularities of $U$.
%
%
\subsection{Multiplier ideals and the openness conjecture}
If $U$ is a quasi-psh function on $X$, then the \emph{multiplier ideal sheaf} $\cJ(U)\subset\cO_X$ is a coherent ideal sheaf whose stalk at any point $x\in X$ is given by the set of holomorphic germs $f\in\cO_{X,x}$ such that $|f|e^{-U}\in\Ltwoloc(X,x)$.\footnote{In the literature, the multiplier ideal sheaf $\cJ(U)$ is usually defined for a psh (rather than quasi-psh) function $U$, but there is no real difference since adding a bounded function $U$ does not affect the multiplier ideal. The same remark applies to many other constructions in this appendix.}
  It was formally introduced by Nadel~\cite{Nadel}. The coherence of $\cJ(U)$ follows from the strong Noetherian property, Krull's Lemma, and H\"ormander's $L^2$-estimates, see the proof of~\cite[Lemma~4.4]{Dem93}.
For general information on multiplier ideals, see~\eg~\cite{DK01} and the references therein.

If $0<\e<\e'$, then $\cJ(U)\supset\cJ((1+\e)U)\supset \cJ((1+\e')U)$, so by the strong Noetherian property of coherent ideals, $\cJ^+(U):=\bigcup_{\e>0}\cJ((1+\e)U)$ is a coherent ideal sheaf on $X$.

The following result is known as the `strong openness conjecture' of Demailly--Koll\'ar~\cite{DK01}, and was established by Guan--Zhou~\cite{GuZh}; see also~\cite{Bernop} for the case $f=1$ and~\cite{valmul,JM14} for the two-dimensional case.
\begin{thm}\label{thm:openness}
  For each quasi-psh function $U$ on $X$, we have $\cJ_+(U)=\cJ(U)$. Equivalently, for any compact subset $K\subset X$, and every holomorphic function $f$ in a neighborhood of $K$ such that $|f|e^{-U}\in\Ltwoloc(X,K)$, there exists $\e>0$ such that $|f|e^{-(1+\e)U}\in\Ltwoloc(X,K)$.  
\end{thm}
The proofs in~\cite{Bernop,GuZh} are based on the Ohsawa--Takegoshi theorem.
%
%
\subsection{More general integrability}\label{S307}
Generalizing the study of multiplier ideals, given a compact subset $K\subset X$ and  quasi-psh functions $U$ and $V$ on a neighborhood of $K$ in $X$, we shall analyze whether $e^{U-V}\in\Ltwoloc(X,K)$.
When studying this problem, we will use the change of variables formula in the following form: 
if $\mu\colon X'\to X$ is a modification, then $e^{U-V}\in\Ltwoloc(X,K)$ iff $e^{U\circ\mu+W_\mu-V\circ\mu}\in\Ltwoloc(X',\mu^{-1}(K))$, where $W_\mu$ is a quasi-psh function on $X'$ with divisorial singularities along $K_{X'/X}$.

Together with Hironaka's theorem in the form above, this allows us, in principle, to determine integrability of $e^{U-V}$ in the case when $U$ and $V$ have analytic singularities. Indeed, passing to a log resolution of the singularities of $U+V$, we reduce to the case when $U$ and $V$ have have divisorial singularities along a common snc divisor. Thus we can cover  $K$ by finitely many coordinate charts $(z_1,\dots,z_n)$ on which $U=\sum_{i=1}^nc_i\log|z_i|+O(1)$ and $V=\sum_{i=1}^nd_i\log|z_i|+O(1)$, where $c_i,d_i\ge 0$, and then $e^{U-V}\in\Ltwoloc(X,K)$ iff $c_i+1>d_i$ for $1\le i\le n$ and each chart.

Calculations similar to the one just performed will appear repeatedly in what follows. For example, they allow us to prove the following generalization of Theorem~\ref{thm:openness}.
\begin{cor}\label{cor:openness}
  Let $K\subset X$ be a compact set, $U,V$ quasi-psh functions defined in a neighborhood of $K$ in $X$, and assume that $U$ has analytic singularities in a neighborhood of $K$.
Then $e^{U-V}\in\Ltwoloc(X,K)$ iff $e^{U-(1+\e)V}\in\Ltwoloc(X,K)$ for some $\e>0$.
\end{cor}
\begin{proof}
  We only need to prove the direct implication as $V$ is bounded above near $K$. After replacing $X$ by a neighborhood of $K$ and applying a log resolution of the singularities of $U$, we reduce to the case when $U$ has divisorial singularities along an snc divisor. By compactness, we may also assume that $K=\{x\}$ is a singleton.
 
  Pick local coordinates $(z_1,\dots,z_n)$ at $x$ such that $U=\sum_{i=1}^nc_i\log|z_i|+O(1)$, where $c_i\ge 0$. If $m$ is an integer with $m\ge\max_ic_i$, then $V':=V+\sum_i(m-c_i)\log|z_i|$ is quasi-psh near $x$, the function $f=\prod_{i=1}^nz_i^m$ is holomorphic near $x$, and $|f|e^{-V'}=e^{U-V}\in\Ltwoloc(X,x)$.
By Theorem~\ref{thm:openness}, there exists $\e>0$ such that $|f|e^{-(1+\e)V'}=e^{U-(1+\e)V-\e(V'-V)}\in\Ltwoloc(X,x)$. This implies $e^{U-(1+\e)V}\in\Ltwoloc(X,x)$ since $V'-V$ is bounded above near $x$.
\end{proof}
%
%
\subsection{Divisorial valuations}
By a \emph{prime divisor over $X$} we mean a connected smooth hypersurface $E\subset X'$, where $X'$ is a complex manifold and $\mu\colon X'\to X$ a modification. The \emph{center} of $E$ on $X$ is defined as $c_X(E):=\mu(E)$.

Two prime divisors $E_1\subset X'_1$ and $E_2\subset X'_2$ over $X$ are \emph{equivalent} if there exist modifications $X''\to X'_i$, $i=1,2$ such that the two compositions $X''\to X'_i\to X$ coincide, and the strict transforms of $E_1$ and $E_2$ on $X''$ is a common smooth hypersurface. By Hironaka's theorem, this defines an equivalence relation. Furthermore, any two non-equivalent prime divisors over $X$ are equivalent, respectively, to two disjoint smooth hypersurfaces in a single $X'$, where $X'\to X$ is a modification.
Equivalent prime divisors have the same center on $X$.

Any prime divisor $E$ over $X$ defines a valuation $\ord_E$ on the field of meromorphic functions on $X$: if $\mu\colon X'\to X$ is a modification, $E\subset X'$ a prime divisor, and $f\not\equiv0$ a meromorphic function on $X$, then $\ord_E(f)\in\Z$ is the order of vanishing of $f\circ\mu$ along $E$. We call $\ord_E$ a \emph{divisorial valuation}.
We can also define $\ord_E(\fa)\in\Z$ for any coherent ideal sheaf $\fa$ on $X$.
The \emph{log discrepancy} of $E$ is $A_X(E):=1+\ord_E(K_{X'/X})$.

Equivalent prime divisors over $X$ induce the same divisorial valuation. The converse is not true in general, since $X$ may not admit any nonconstant meromorphic functions. Similarly, if $\tX$ is an open subset of $X$, a prime divisor over $\tX$ may not extend to a prime divisor over $X$, and a prime divisor $E$ over $X$ with $c_X(E)\cap\tX\ne\emptyset$ does not define a unique prime divisor over $\tX$, in general.

We have a better understanding of prime divisors over $X$ whose center on $X$ is a \emph{point} $x\in X$. We will use the following piece of (non-standard) terminology. A  \emph{blowup over $x$} is the blowup $X'\to X$ of a coherent ideal sheaf $\fa\subset\cO_X$ co-supported at $x$, such that $X'$ is a complex manifold. This is the same thing as a projective modification $\mu\colon X'\to X$ (with $X'$ a complex manifold) that is an isomorphism above $X\setminus\{x\}$. We also note that coherent ideal sheaves on $X$ co-supported at $x$ are in 1-1 correspondence with $\fm_x$-primary ideals of $\cO_{X,x}$, where $\fm_x$ is the maximal ideal of $\cO_{X,x}$.
\begin{lem}\label{L305}
  If $x\in X$ is a point, then:
  \begin{itemize}
  \item[(i)]
    any prime divisor $E$ over $X$ with $c_X(E)=x$ induces a valuation $\ord_E\colon\cO_{X,x}\to\Z_{\ge0}$ that is strictly positive on the maximal ideal $\fm_x\subset\cO_{X,x}$;
  \item[(ii)]
    two prime divisors $E_1,E_2$ over $X$ with $c_X(E_1)=c_X(E_2)=x$ are equivalent iff they induce the same valuation on $\cO_{X,x}$;
  \item[(iii)]
    any prime divisor over $X$ with center $x$ is equivalent to a connected smooth hypersurface on a blowup over $x$.
  \end{itemize}
\end{lem}
\begin{proof}
  The statement in~(i) is clear, as is the direct implication in~(ii).

  Let $\mathrm{EPDiv}(X,x)$ be the set of equivalence classes of prime divisors over $X$ with center $x$. For any open neighborhood $U$ of $x$ in $X$, we have a canonical injective map $\mathrm{EPDiv}(X,x)\to\mathrm{EPDiv}(U,x)$. Also note that by Hironaka's theorem, for any modification $\mu\colon X'\to X$, there exists an open neighborhood $U$ of $x$ and a projective modification $U''\to U$ such that the induced bimeromorphic map $U''\to \mu^{-1}(U)$ is holomorphic.

  To prove the reverse implication in~(ii) we may therefore assume that there exists a projective modification $\mu\colon X'\to X$ with $E_1,E_2$ distinct (even disjoint) smooth hypersurfaces on $X'$. After shrinking $X$, if necessary, there exists a $\mu$-ample divisor $A$ on $X$. Pick $m\gg1$ such that the line bundles $\cO_{X'}(mA)$ and $\cO_{X'}(mA-E_1)$ are both $\mu$-globally generated. The quotient of two general sections of these line bundles defines an element $f$ in the fraction field of $\cO_{X,x}$ such that $\ord_{E_1}(f)\ne\ord_{E_2}(f)$. This implies $\ord_{E_1}\ne\ord_{E_2}$ on $\cO_{X,x}$.
  
  It remains to prove~(iii). If $U$ is an open neighborhood of $x$ in $X$, then any blowup of $U$ over $x$ induces a unique blowup of $X$ over $x$. Together with the observations above, we only need to consider the case of a prime divisor $E$ over $X$ associated to a projective modification. The induced valuation $\ord_E$ on $\cO_{X,x}$ is divisorial: its value group is $\Z$ and it has transcendence degree $n-1$. It therefore follows from~\cite[Proposition~3.7]{JM12} that there exists a projective birational morphism (of schemes) $\mu\colon V\to\Spec\cO_{X,x}$ that is an isomorphism over the complement of the closed point, and a prime divisor $F$ on $V$ such that $\ord_E$ is the order of vanishing along $F$. We may assume that $F$ is smooth. Since $V$ is smooth, $\mu$ is the blowup along an $\fm_x$-primary ideal $\fa\subset\cO_{X,x}$, which we can view as coherent ideal sheaf on $X$. By blowing up $X$ along $\fa$, we obtain a blowup $Y'\to X$ over $x$ and a connected hypersurface $F'\subset Y$ (corresponding to $F$) such that $\ord_{F'}=\ord_E$ on $\cO_{X,x}$. By~(ii) it follows that $E$ and $F'$ are equivalent.
\end{proof}

%
%
\subsection{Divisorial valuations and Lelong numbers}
Let $E$ be a prime divisor over $X$.
Following~\cite{hiro} we can define $\ord_E(U)\in\R_{\ge0}$ for any quasi-psh function $U$ on $X$. Suppose $\mu\colon X'\to X$ is a modification such that $E\subset X'$ a connected smooth hypersurface. Then $U\circ\mu$ is a quasi-psh function on $X'$. Let $U_E$ be a quasi-psh function on $X'$ with divisorial singularities of type $E$. We define 
$\ord_E(U)$ as the \emph{generic Lelong number} of $U$ along $E$, that is, the supremum of all numbers $t\ge0$ such that $U\circ\mu\le tU_E+O(1)$ near any point $x'\in E$. In fact, we then have $U\circ\mu\le \ord_E(U)U_E+O(1)$ near $E$, and even a \emph{Siu decomposition} $U\circ\mu=\ord_E(U)U_E+U'$, where $U'$ is quasi-psh.

Note that $\ord_E(U)$ only depends on the equivalence class of $E$ Also note that 
if $U$ has analytic singularities of type $\prod_i\fa_i^{c_i}$, then $\ord_E(U)=\sum_ic_i\ord_E(\fa_i)$.
When $E$ is the exceptional divisor of the blowup of a point $x\in X$, $\ord_E(U)$ is the usual Lelong number of $U$ at $x$. In general, $\ord_E(U)$ equals the Lelong number of $U\circ\mu$ at a very general point on $E$.

\begin{lem}\label{lem:lesssing}
  Let $U$ be a quasi-psh function on $X$. Then
  \begin{equation}\label{e303}
    \ord_E(\cJ(U))\le\ord_E(U)<\ord_E(\cJ(U))+A_X(E)
  \end{equation}
  for any prime divisor $E$ over $X$.
\end{lem}
\begin{proof}
  The first inequality follows from the Ohsawa--Takegoshi theorem as in the proof of a celebrated regularization theorem due to Demailly. Namely, it follows from the proof of~\cite[Theorem~4.2~(3)]{DK01} that for any $x\in X$  we can pick generators $f_1,\dots,f_m$ of $\cJ(U)_x$ such that $U\le\max_j\log|f_j|+O(1)$ near $x$. This implies the first inequality of the proposition.

  The second inequality follows from a direct computation as in~\S\ref{S307}. Let $V$ be a quasi-psh function with analytic singularities of type $\cJ(U)$. Thus $e^{V-U}\in\Ltwoloc(X)$. Pick a modification $\mu\colon X'\to X$ such that $E\subset X'$ is a smooth hypersurface and $V\circ\mu$ has divisorial singularities along an snc divisor whose support includes $E+K_{X'/X}$. Let $W_\mu$ be a quasi-psh function with divisorial singularities along $K_{X'/X}$, and set $V'=V\circ\mu+W_\mu$, $U'=U\circ\mu$. The change of variables formula shows that $e^{V'-U'}\in\Ltwoloc(X)$. Let $(z_1,\dots,z_n)$ be local coordinates at a general point $x'\in E$ such that $E=\{z_1=0\}$. Then 
  $U'\le\ord_E(U)\log|z_1|+O(1)$ and $V'=(\ord_E(V)+A_X(E)-1)\log|z_1|+O(1)$ near $x'$, and the desired inequality follows.
\end{proof}
%
%
\subsection{A valuative criterion of integrability}
The following theorem is the main result in this appendix. It yields a valuative description of the multiplier ideal of a general quasi-psh function.
\begin{thm}\label{thm:valcritint1}
  Let $X$ be a complex manifold, $K\subset X$ a compact subset, and $U,V$ quasi-psh functions defined in a neighborhood of $K$ on $X$. Assume that $U$ has analytic singularities in a neighborhood of $K$. Then the following assertions are equivalent:
  \begin{itemize}
  \item[(i)]
    $e^{U-V}\in\Ltwoloc(X,K)$;
  \item[(ii)]
    there exists $\e>0$ such that 
    \begin{equation}\label{e302}
      \ord_E(U)+A_X(E)\ge(1+\e)\ord_E(V)
    \end{equation}
    for all prime divisors $E$ over $X$ with $c_X(E)\subset K$;
  \item[(iii)]
    there exists $\e>0$ such that~\eqref{e302} holds for all prime divisors $E$ over $X$ such that  $c_X(E)$ is a point contained in $K$.
  \end{itemize}
\end{thm}
\begin{rmk}\label{R306}
  Condition~(i) is unchanged by replacing $X$ by an neighborhood of $K$ in $X$. The same is therefore true for~(ii) and~(iii). In the case of~(iii), this can also be seen from Lemma~\ref{L305}.
\end{rmk}
\begin{proof}
  Let $\tX$ be a neighborhood of $K$ in $X$ such that $U$ and $V$ are defined on $\tX$, and $U$ has analytic singularities there. It suffices to prove the statement with $X$ replaced by $\tX$. Indeed,~(i) is not affected by passing to $\tX$, and the same is true for~(iii) in view of Lemma~\ref{L305}. As for~(ii), the version on $\tX$ is evidently stronger than the one on $X$ (any prime divisor $E$ over $X$ with $c_X(E)\subset K$ induces a unique prime divisor over $\tX$) and both are stronger than~(iii).
  Thus we may assume that $U$ and $V$ are defined on $X$, and $U$ has analytic singularities there.

  Next we show that~(i)--(iii) are invariant under modifications. 
  Consider a modification $\mu\colon X'\to X$. Let $W_\mu$ be a quasi-psh function on $X'$ with divisorial singularities of type $K_{X'/X}$, and set $U':=U\circ\mu+W_\mu$ and $V':=V\circ\mu$. Then $U'$ has analytic singularities, and by the change of variables formula we see that~(i) is equivalent to $e^{U'-V'}\in\Ltwoloc(X,K')$, where $K'=\mu^{-1}(K)$.
  We can also identify the set of (equivalence classes of) prime divisors $E$ over $X$ and $X'$, and we have $\ord_E(U)+A_X(E)=\ord_E(U')+A_{X'}(E)$, $\ord_E(V)=\ord_E(V')$, and $c_X(E)=\mu(c_{X'}(E))$. Thus we may assume that $U$ has divisorial singularities along an snc divisor.
  
  \smallskip
  We now introduce an auxiliary condition:
  \begin{itemize}
  \item[(iv)]
    there exists $\e>0$ such that~\eqref{e302} holds for all prime divisors $E$ over $X$ with $c_X(E)\cap K\ne\emptyset$.
  \end{itemize}
  
  \smallskip
  Evidently, (iv)$\Rightarrow$(ii)$\Rightarrow$(iii), so it suffices to prove
  (iii)$\Rightarrow$(iv), (i)$\Rightarrow$(ii), and (iv)$\Rightarrow$(i).
  
  \smallskip
  To prove that (iii) implies~(iv) we use an approximation procedure. Consider a prime divisor $E$ over $X$ with $c_X(E)\cap K\ne\emptyset$. Passing to a modification as above, we may assume $E\subset X$ is a connected smooth hypersurface.
  Consider a point $x\in c_X(E)\cap K$, and pick local coordinates $(z_1,\dots,z_n)$ centered at $x$ such that $E=\{z_1=0\}$. For any $m\ge 1$, consider the monomial valuation $w_m$ on $\cO_{X,x}$ in these coordinates with $w_m(z_1)=1$ and $w_m(z_j)=1/m$ for $2\le j\le n$. There exists a divisor $E_m$ over $X$ such that $\ord_{E_m}=mw_m$. Indeed, $E_1$ is the exceptional divisor of the blowup of $x$, and, for $m\ge 2$, $E_m$ is the exceptional divisor of the blowup of the intersection of $E_{m-1}$ and the strict transform of $E$.

  Now $c_X(E_m)=x$ for all $m$, and $w_m$ is an approximation of $w=\ord_E$ in the following sense. First, $A_X(w_m)=1+(n-1)/m$ and $A_X(w)=1$, so $\lim_{m\to\infty}A_X(w_m)=A_X(w)$.
  Second, for any ideal $\fb\subset\cO_{X,x}$ we have $\lim_{m\to\infty}w_m(\fb)=w(\fb)$. This immediately shows that $w_m(U)\to w(U)$ as $m\to\infty$ since $U$ has analytic singularities. Finally, it follows from~\eqref{e303} and what precedes that $\lim_{m\to\infty}w_m(V)=w(V)$. Since~\eqref{e302} holds for $E_m$, it must also hold for $E$.
  
  \smallskip
  To prove~(i)$\Rightarrow$(ii), we use the openness conjecture in the form of Corollary~\ref{cor:openness}. Thus assume $e^{U-(1+\e)V}\in\Ltwoloc(X,K)$, where $\e>0$. We must prove~\eqref{e302} for any prime divisor $E$ over $X$ such that $c_X(E)\subset K$. Passing to a modification, we may assume that $E$ is a connected smooth hypersurface on $X$ with $E\subset K$, and that $U$ has divisorial singularities along an snc divisor whose support contains $E$. Then we can argue as in the proof of~\eqref{e303}.
  Pick a general point $x\in E$ such that $U$ and $V$ are defined in a neighborhood of $x$, and let $(z_1,\dots,z_n)$ be local analytic coordinates at $x$ such that $E=\{z_1=0\}$. Then $U=\ord_E(U)\log|z_1|+O(1)$ and $V\le\ord_E(V)\log|z_1|+O(1)$ near $x$. Since $e^{U-(1+\e)V}\in\Ltwoloc(x)$ and $A_X(E)=1$,~\eqref{e302} must hold (with strict inequality).

  \smallskip
  It remains to prove (iv)$\Rightarrow$(i). First suppose that $U$ and $V$ both have analytic singularities. After passing to a suitable modification, we may assume $U$ and $V$ both have divisorial singularities along a common reduced snc divisor $D$. We again follow the proof of~\eqref{e303}. Pick any point $x\in K$, and let $E_1,\dots,E_m$, where $0\le m\le n$, be the irreducible components of $D$ that contain $x$. Pick local coordinates $(z_1,\dots,z_n)$ at $x$ such that $E_i=\{z_i=0\}$ for $1\le i\le n$. Then $U=\sum_{i=1}^mc_i^{(U)}\log|z_i|+O(1)$ and 
  $V=\sum_{i=1}^mc_i^{(V)}\log|z_i|+O(1)$ near $x$, where $c_i^{(U)},c_i^{(V)}\ge0$. Now~\eqref{e302} gives $(1+\e)c_i^{(V)}\le c_i^{(U)}+1$ for all $i$. Since $c_i^{(U)},c_i^{(V)}\ge0$, this implies $c_i^{(V)}-c_i^{(U)}>-1$ for all $i$, and hence $e^{U-V}\in\Ltwoloc(X,x)$.

  Now consider the general case. For each $k\in\Z_{>0}$, let $V_k$ be a quasi-psh function on $X$ with analytic singularities of type $\cJ(kV)^{1/k}$. Combining~(ii) and Lemma~\ref{lem:lesssing}, we have 
\begin{equation*}
  \ord_E(U)+A_X(E)\ge(1+\e)\ord_E(V)\ge(1+\e)\ord_E(V_k)
\end{equation*}
for all $k$ and all prime divisors $E$ over $X$ with $c_X(E)\cap K\ne\emptyset$.
By what precedes, this implies $e^{U-(1+\e)V_k}\in\Ltwoloc(X,K)$.
Pick $k\ge1$ such that $\frac1{k-1}\le\e$. Then $e^{\frac{k}{k-1}(U-V_k)}\in\Ltwoloc(X,K)$, since $U$ and $V_k$ are bounded above near $K$. On the other hand, the definition of $V_k$ gives $e^{k(V_k-V)}\in\Ltwoloc(X,K)$, and we conclude using H\"older's inequality.
\end{proof}
%
%
\subsection{The algebraic case}
Now assume that $X$ is a complex algebraic manifold, \ie the analytification of a smooth complex algebraic variety. We will translate the integrability criterion in Theorem~\ref{thm:valcritint1} to one using the set $\Xdiv$ of (rational) divisorial valuations on $X$. As we need it in the proof of Theorem~\ref{thm:Lasymp}, we do this in the presence of a (possibly trivial) $\C^*$-action on $X$.

Any (nontrivial) valuation $w\in\Xdiv$ is of the following form: there exists a projective birational morphism $X'\to X$, with $X'$ a complex algebraic manifold, and a smooth hypersurface $E\subset X'$ such that $w=c\ord_E$ for some $c\in\Q_{>0}$. Note that $E$ defines a prime divisor over $X$ in the sense above. 
As a partial converse, if follows from Lemma~\ref{L305} that any prime divisor over $X$ whose center on $X$ is a (closed) point is equivalent to a divisor over $X$ obtained from an element of $\Xdiv$.

\begin{thm}\label{thm:valcritint2}
  Let $X$ be a complex algebraic manifold with a $\C^*$-action, and $Y\subset X$ a compact Zariski closed subset that is invariant under this action. Let $U,V$ be $S^1$-invariant quasi-psh functions defined in a neighborhood of $Y$ on $X$. Assume that $U$ has analytic singularities in a neighborhood of $Y$. Then the following assertions are equivalent:
  \begin{itemize}
  \item[(i)]
    $e^{U-V}\in\Ltwoloc(X,Y)$;
  \item[(ii)]
    there exists $\e>0$ such that
    \begin{equation}\label{e306}
      w(U)+A_X(w)\ge(1+\e)w(V)
    \end{equation}
    for all $w\in\Xdiv$ with $c_X(w)\subset Y$;
  \item[(iii)]
    there exists $\e>0$ such that~\eqref{e306} holds for all $\C^*$-invariant valuations $w\in\Xdiv$ such that  $c_X(w)$ is a point contained in $Y$.
  \end{itemize}
\end{thm}
\begin{proof}[Proof of Theorem~\ref{thm:valcritint2}]
  The implication~(ii)$\Rightarrow$(iii) is clear, and (i)$\Rightarrow$(ii) is a special case of Theorem~\ref{thm:valcritint1}. It remains to prove (iii)$\Rightarrow$(i). To do so, let $\tX$ be an $S^1$-invariant open neighborhood of $Y$ in $X$ on which $U$ and $V$ are defined, and $U$ has analytic singularities. Consider the following condition:
  \begin{itemize}
  \item[(iv)]
    there exists $\e>0$ such that $\ord_E(U)+A_{\tX}(E)\ge(1+\e)\ord_E(V)$ for all $S^1$-invariant prime divisors $E$ over $\tX$ such that $c_{\tX}(E)\cap Y\ne\emptyset$.
  \end{itemize}
  By the condition on $E$ we mean the following: there exists a modification $\mu\colon X'\to\tX$ such that the $S^1$-action on $\tX$ lifts to $X'$, and $E$ is a connected smooth $S^1$-invariant hypersurface of $X'$.

  We can then prove~(iv)$\Rightarrow$(i) in the same way as~(ii)$\Rightarrow$(i) in the proof of Theorem~\ref{thm:valcritint1}. Indeed, the $S^1$-invariance of $U$ and $V$ implies that the functions $V_k$ can also be taken $S^1$-invariant. To test the integrability of $e^{2(U-\frac{k}{k-1}V_k)}$ in a neighborhood of $Y$, it suffices to use $S^1$-invariant prime divisors over $X$, as we can find a common $S^1$-equivariant log resolution of the singularities of $U$ and $V_k$.

  We finally prove that~(iii) implies~(iv), using a variant of the argument in the proof of~(iii)$\Rightarrow$(ii) in Theorem~\ref{thm:valcritint1}. Let $\tX$ be an $S^1$-invariant neighborhood of $Y$ in $X$, $\mu\colon X'\to\tX$ a modification such that the $S^1$-action lifts to $X'$, and $E$ is a connected smooth $S^1$-invariant hypersurface of $X'$ such that $\mu(E)\cap Y\ne\emptyset$.
  Let $Z$ be an irreducible component of $E\cap\mu^{-1}(Y)$. Then $Z\subset X$ is a proper subvariety that is $S^1$-invariant, and hence $\C^*$-invariant. By the valuative criterion of properness, there exists a (closed) point $x'\in E\cap\mu^{-1}(Y)$ that is fixed under the $\C^*$-action.
  
  Pick local coordinates $(z_1,\dots,z_n)$ at $x'$ such that $E=\{z_1=0\}$. Since $E$ and $x'$ are $S^1$-invariant, the monomial valuation $w_m$ with weights $w_m(z_1)=1$ and $w_m(z_j)=1/m$, $2\le j\le n$ is also $S^1$-invariant, and hence $\C^*$-invariant. Note that $w_m$ is an element of $\Xdiv$ in view of Lemma~\ref{L305}.
  The inequality in~\eqref{e306} for $w_m$ now implies the inequality in~(iv) as $m\to\infty$.
\end{proof}

%
%
\subsection{Twisted log-canonical thresholds}
Now let $X$ be a smooth projective variety and $\theta$ a quasi-positive $(1,1)$-current on $X$, that is, $\theta=\theta_0+dd^c\p$, where $\theta_0$ is smooth and $\p$ quasi-psh on $X$. We then define $\cJ(\theta):=\cJ(\p)$. By definition, $\theta$ is klt iff $\cJ(\theta)=\cO_X$.

For $v\in\Xdiv$, we set
\begin{equation*}
  A_\theta(v)
  :=A_X(v)-v(\theta),
\end{equation*}
where $v(\theta):=v(\p)$. (This does not depend on the decomposition $\theta=\theta_0+dd^c\p$.)
\begin{cor}\label{cor:lctval1}
  The current $\theta$ is klt iff there exists $\e>0$ such that $A_\theta\ge\e A_X$ on $\Xdiv$.
\end{cor}  
\begin{proof}
  In the notation above, $\theta$ is klt iff $e^{-2\p}$ is locally integrable at any point on $x$. It follows from Theorem~\ref{thm:valcritint1} that this is the case iff there exists $\e'>0$ such that $0\ge(1+\e')v(\p)-A_X(v)$ for every $v\in\Xdiv$. We can then take $\e=\e'/(1+\e')$.
\end{proof}
Now suppose $\theta$ is a klt current. 
If $D$ is an effective $\Q$-divisor on $X$, then we write $D$ also for the current of integration of $D$; this allows us to define $\cJ(\theta+D)$. The \emph{log canonical threshold} of $D$ with respect to $\theta$ is now defined as
\begin{equation*}
  \lct_\theta(D):=\sup\{c\in\Q_{\ge0}\mid \cJ(\theta+cD)=\cO_X\}.
\end{equation*}
\begin{cor}\label{cor:lctval2}
  With notation as above, we have
  \begin{equation}\label{e304}
    \lct_\theta(D)=\inf_{v\in\Xdiv\setminus\{v_\triv\}}\frac{A_\theta(v)}{v(D)},
  \end{equation}
  where $A_\theta(v)=A_X(v)-v(\theta)$.
\end{cor}
\begin{proof}
  We may assume $D\ne0$ or else the equality holds with both sides equal to $+\infty$.
  Pick $\e_0>0$ such that $\e_0v(\theta)\le A_\theta(v)$ for all $v\in\Xdiv$; this is possible in view of Corollary~\ref{cor:lctval1}.
  
  If $c<\lct_\theta(D)$, then Theorem~\ref{thm:valcritint1} shows that there exists $\e>0$ such that $0\ge (1+\e)(v(\theta)+cv(D))-A_X(v)$, and hence $A_\theta(v)\ge\e v(\theta)+(1+\e)cv(D)\ge cv(D)$ for all $x\in\Xdiv$.

  Conversely, suppose $c\in\Q_{>0}$ and $c>\lct_\theta(D)$, so that $\cJ(\theta+cD)\ne\cO_X$. Pick $\e'\in\Q_{>0}$ such that $(1-\e'/\e_0)\ge c(1+\e')$. By Theorem~\ref{thm:valcritint1} we can find $v\in\Xdiv$ such that $A_X(v)<(1+\e')(v(\theta)+cv(D))$. Together with the inequality $\e_0v(\theta)\le A_\theta(v)$ and the choice of $\e'$, this implies $A_\theta(v)\le cv(D)$, which completes the proof.
\end{proof}
%
%
%
%

\end{document}